\documentclass[dvipsnames]{m2an}
\usepackage{mypackages}
\usepackage{mycommands}
\usepackage{plot_options}
\begin{document}
%%-----------------------------
%%      the top matter
%%-----------------------------
\title{Explicit stabilized multirate methods for the monodomain model in cardiac electrophysiology}
\thanks{This work was supported by the European High-Performance Computing Joint Undertaking EuroHPC under grant agreements No 955495 (MICROCARD) and No 955701 (TIME-X), co-funded by the Horizon 2020 programme of the European Union (EU) and the Swiss State Secretariat for Education, Research and Innovation. We also acknowledge the CSCS-Swiss National Supercomputing Centre (project no.~s1074)}
%\thanks{...}% At most 5 thankste
%
\author{Giacomo Rosilho de Souza}
\address{Euler Institute, Università della Svizzera italiana, via G.~Buffi 13, 6900 Lugano, Switzerland;
\\ \email{\{giacomo.rosilhodesouza,rolf.krause\}@usi.ch}}
\author{Marcus J.~Grote}
\address{Department of Mathematics, University of Basel, Rheinsprung 21, 4051 Basel, Switzerland;\\
\email{marcus.grote@unibas.ch}}
\author{Simone Pezzuto}
\sameaddress{1}\secondaddress{Laboratory for Mathematics in Biology and Medicine, Department of Mathematics, Università di Trento, via Sommarive 14, 38123, Trento, Italy; \email{simone.pezzuto@unitn.it}}
\author{Rolf Krause}
\sameaddress{1}\secondaddress{Faculty of Mathematics and Informatics, FernUni Schweiz, Schinerstrasse 18, 3900, Brig, Switzerland}
\runningauthors{G. Rosilho de Souza \etal}
\date{dates will be set by the publisher}

\begin{abstract}
Fully explicit stabilized multirate (mRKC) methods are well-suited for the numerical solution of large multiscale systems of stiff ordinary differential equations thanks to their improved stability properties. To demonstrate their efficiency for the numerical solution of stiff, multiscale, nonlinear parabolic PDE's, we apply mRKC methods to the monodomain equation from cardiac electrophysiology. In doing so, we propose an improved version, specifically tailored to the monodomain model, which leads to the explicit exponential multirate stabilized (emRKC) method. Several numerical experiments are conducted to evaluate the efficiency of both mRKC and emRKC, while taking into account different finite element meshes (structured and unstructured) and realistic ionic models. The new emRKC method typically outperforms a standard implicit-explicit baseline method for cardiac electrophysiology. Code profiling and strong scalability results further demonstrate that emRKC is faster and inherently parallel without sacrificing accuracy.
\end{abstract}
 \begin{resume} 
% French version.
Les m\'ethodes multirate stabilis\'ees explicites (mRKC) sont bien adapt\'ees à la r\'esolution num\'erique de grands syst\`emes multi-\'echelles d'\'equations diff\'erentielles ordinaires rigides gr\^ace \`a leur meilleure stabilit\' e. Afin de démontrer leur efficacité pour la résolution numérique d'EDP paraboliques non-linéaires raides et multi-échelles, nous appliquons les méthodes mRKC à l'équation monodomaine de l'électro\-physio\-logie cardiaque. 
Ce faisant, nous proposons une version améliorée, spécifiquement adaptée au modèle monodomaine, qui conduit à la m\'ethode exponentielle multirate stabilisée (emRKC) explicite.
Plusieurs expériences numériques permettent d'évaluer l'efficacité des m\'ethodes mRKC et emRKC, tout en tenant compte de différents maillages d'éléments finis (structurés et non-structurés) 
et des modèles ioniques réalistes. La nouvelle méthode emRKC
surpasse généralement une méthode de base implicite-explicite standard pour l'électrophysiologie cardiaque.
Le profilage du code et les résultats de scalabilité démontrent en outre que emRKC est plus rapide et intrinsèquement parallèle sans sacrifier la précision.
 \end{resume}

\maketitle

% \GR{Marcus, you can add comments with the \texttt{\textbackslash MG\{your comment\}} command.}
% \GR{Simone, with the \texttt{\textbackslash SP\{your comment\}} command.}

\section{Introduction}\label{sec:intro}

Spatial discretizations of parabolic partial differential equations (PDEs), possibly nonlinear,  typically lead to
large systems of stiff ordinary differential equations (ODEs). The degree of stiffness, be it mild or severe, depends in particular on the rate of change of individual components or the mesh size. For time discretization, there is generally a trade-off between explicit and implicit time-marching schemes. Explicit schemes are 
easier to implement and also cheaper per time-step, but require more steps due to stability constraints. In contrast, implicit schemes require sophisticated (parallel) solvers for the needed linear
or nonlinear systems and thus are more expensive per time-step, but they also require fewer time-steps. Choosing between explicit and implicit time integration is a delicate process which depends on several factors, such as the stiffness of the system or the efficiency of the linear or nonlinear solvers available. Indeed,
an unfortunate choice can increase the computational cost by orders of magnitude due to abrupt changes in any methods' performance.

Explicit stabilized methods, such as Runge--Kutta--Chebyshev (RKC) methods \cite{HoS80,SSV98}, 
strike a balance between explicit and implicit schemes while their computational cost smoothly depends on the degree of stiffness. On the one hand, they are fully explicit and thus simple to implement in parallel, even in the presence of nonlinearity; on the other hand, they adapt to the current degree of stiffness by modifying the number of stages accordingly. Since their stability domain grows quadratically with the number of stages, they are significantly more efficient than standard explicit methods.
%; in fact, their computational cost is only the square root of the cost of the former.
Compared to implicit methods, they can also be more efficient in particular for large-scale three-dimensional linear, or nonlinear, problems\cite{AbM01,Abd02,AbV13,DDD13}. In summary, explicit stabilized schemes combine the advantages of both
explicit and implicit methods while avoiding any abrupt performance deterioration when transitioning from a non-stiff to stiff regime, or from linear to nonlinear problems.

When stiffness is induced only by a few components, as in the presence of spatial local mesh
refinement or coupled chemical reactions, the efficiency of RKC methods can deteriorate. 
To overcome that crippling effect due to a few severely stiff components,
multirate Runge–Kutta–Chebyshev (mRKC) methods were proposed in \cite{AGR22};
they remain fully explicit, and thus easy to parallelize, 
though with a stability condition independent of those few severely stiff components.
Explicit multirate RKC methods thus permit to circumvent any overly stringent 
stability constraint due to local mesh refinement; hence, they also
extend local time-stepping (LTS) methods \cite{grote2009,Grote2010,GroteMitkova13,GMM15} to parabolic PDE's.

Here we shall consider the monodomain model for cardiac electrophysiology, which consists of a parabolic PDE for the diffusion of the electric potential and a system of ODE's for the ionic model of the intracellular and membrane dynamics of the cardiac cell~\cite{colli2014mathematical}. 
After spatial discretization, the monodomain model can be written as
\begin{equation}\label{eq:FSEeq}
	y'=f_F(y)+f_S(y)+f_E(y),  \qquad y(0)=y_0,
\end{equation}
where $f_F$ is cheap but stiff (``F'' for ``fast''), $f_S$ is non-stiff but expensive (``S'' for ``slow''), and $f_E$ is in general a severely stiff term of the form
\begin{equation}\label{eq:def_fE}
f_E(y)=\Lambda(y)(y-y_\infty(y)),
\end{equation}
with $\Lambda(y)$ a diagonal matrix and $y_\infty(y)$ a nonlinear function.
Typically, $f_F$ represents the discrete diffusion operator while $f_S$, $f_E$ depend on the ionic model; here, $f_E$ corresponds to the gating variables and $f_S$ to the remaining ones. 

A common and efficient approach for the time integration of the non-diffusive component of \cref{eq:FSEeq} is based on the Rush--Larsen scheme~\cite{RuL78,Coudiere2020}. The Rush--Larsen scheme is a first-order splitting method employing explicit Euler for $f_S$ and exponential Euler for $f_E$ (hence the ``E'' notation for the ``exponential'')-- see \cite{LinGerJahLoeWeiWie23} for a recent review paper on the possible variants of this approach.
The diffusive term $f_F$ is treated either explicitly or implicitly, depending on the degree of stiffness. For the monodomain model, $f_F$ is mildly stiff and, as a consequence, the choice between an explicit or an implicit approach somewhat challenging. Clearly, explicit stabilized schemes, such as RKC methods, can efficiently treat the $f_S$ and $f_F$ terms, but often perform poorly for the full monodomain model, due to its multiple components and inherent multiscale nature.

In contrast, multirate explicit stabilized methods such as mRKC \cite{AGR22,CrR22} are more appropriate for the monodomain model because they can adapt the number of stages locally to individual components
depending on their respective degree of stiffness.
%Indeed, they have already been introduced in the context of multiscale ODEs. 
Here following \cite{RuL78,Coudiere2020}, we shall further take full advantage of the special structure \cref{eq:def_fE} of $f_E$ by 
%In this paper, we introduce an adaptation of mRKC that takes full advantage of the special structure %\cref{eq:def_fE} of $f_E$. Inspired by the Rush--Larsen scheme, we 
integrating it with an exponential method, while integrating $f_F$ and $f_S$ explicitly with the mRKC scheme. We denote this combined exponential--mRKC time integrator by \emph{exponential multirate RKC} (emRKC). The emRKC method is cheap because it is explicit while utilizing an exponential integrator for the most severely stiff yet diagonal terms. Albeit only first-order, it is more accurate than other commonly used first-order schemes, as its special structure, inherited from mRKC, guarantees a robust coupling between the different model components. Although the emRKC method was developed specifically for the monodomain equation (for both ventricular or atrial geometries, or even the entire heart), it can be applied to any physical model that can be written as \cref{eq:FSEeq,eq:def_fE}.

In \cref{sec:monodomain}, we first present the monodomain model from cardiac electrophysiology. Next, in \cref{sec:mES}, we recall the RKC \cite{HoS80,SSV98} and mRKC methods from \cite{AGR22,CrR22}, summarize their properties, and describe them from an algorithmic point of view. In \cref{sec:exp_mES}, we first adapt mRKC to the particular structure \cref{eq:FSEeq} of the monodomain model, which yields the new emRKC method; then we analyze its accuracy and stability properties. In \cref{sec:num}, we compare both mRKC and emRKC to a standard baseline method on a sequence of numerical experiments. Both two- and three-dimensional problems are considered, either with structured or unstructured meshes, including realistic ionic models. Concluding remarks are presented in \cref{sec:conclusion}.

\section{Monodomain model}\label{sec:monodomain}
The monodomain model is commonly used to simulate electric potential propagation in cardiac tissue. Though a simplification of the more complex bidomain model, it remains important in practice. In this section, we provide a brief introduction to the model; for a more complete description, we refer to \cite{CBC11}.
%\subsection{Monodomain model for cardiac electrophysiology}\label{sec:monodomain_equations}
The monodomain model is given by the following reaction-diffusion system of equations:
\begin{subequations}\label{eq:monodomain}
\begin{alignat}{2}
C_m \frac{\partial V_m}{\partial t} &= \chi^{-1} \nabla \cdot \bigl(\mathbf{D} \nabla V_m\bigr) + I_\text{stim}(t,\bm x) - I_\text{ion}(V_m,\bm{z}_E,\bm{z}_S), \qquad && \text{in $\Omega\times (0,T]$,} \label{eq:monodomain_pde} \\
\frac{\partial \bm{z}_E}{\partial t} &= \bm{g}_E(V_m,\bm{z}_E) = \bm{\alpha}(V_m)\bigl(1-\bm{z}_E\bigr) - \bm{\beta}(V_m) \bm{z}_E, && \text{in $\Omega\times (0,T]$,} \label{eq:monodomain_gate} \\
\frac{\partial\bm{z}_S}{\partial t} &= \bm{g}_S(V_m,\bm{z}_E,\bm{z}_S), && \text{in $\Omega\times (0,T]$,} \label{eq:monodomain_g} \\
-\mathbf{D}\nabla V_m\cdot \bm{n} &= 0, && \text{on }\partial \Omega\times (0,T], \\
V_m(0,\bm x) &= V_{m,0}, \quad {\bm z}_E(0,\bm x) = {\bm z}_{E,0}, \quad {\bm z}_S(0,\bm x) = {\bm z}_{S,0}, && \text{in $\Omega$},
\end{alignat}
\end{subequations}
where parameters and variables are described in \cref{tab:monopar}.
Electric conduction in the myocardium is anisotropic, with a conductivity tensor of the form
\begin{equation}\label{eq:def_D}
\mathbf{D}(\bm x) = 
  \sigma_\ell \bm{a}_\ell(\bm x) \otimes \bm{a}_\ell(\bm x)
+ \sigma_t    \bm{a}_t(\bm x) \otimes \bm{a}_t(\bm x)
+ \sigma_n    \Bigl( \mathbf{I} - \bm{a}_\ell(\bm x) \otimes \bm{a}_\ell(\bm x) - \bm{a}_t(\bm x) \otimes \bm{a}_t(\bm x) \Bigr),
\end{equation}
where $\sigma_\star = \sigma_\star^i\sigma_\star^e/(\sigma_\star^i+\sigma_\star^e)$ and $\sigma^i_\star$, $\sigma^e_\star$ are also given in \cref{tab:monopar}. The vector fields $\bm a_\ell(\bm x)$ and $\bm a_t(\bm x)$ are orthonormal and are aligned with the local fiber and sheet direction, respectively.
%V_m(t,\bm x)$ (\si[per-mode=symbol]{\milli\volt}) is the transmembrane potential, $\chi(\bm x)$ (\si{\per\milli\metre}) is the surface area-to-volume ratio of cardiac cells, $C_m$ (\si[per-mode=symbol]{\micro\farad\per\milli\metre\squared}) is the cell membrane capacitance, $I_\text{ion}$ (\si[per-mode=symbol]{\micro\ampere\per\milli\metre\squared}) is the transmembrane current density, $I_\text{stim}(t,\bm x)$ (\si[per-mode=symbol]{\micro\ampere\per\milli\metre\squared}) is the stimulus current, $\mathbf{D}(\bm x)$ (\si[per-mode=symbol]{\milli\siemens\per\milli\metre}) is the electric conductivity tensor, and  $\bm{n}$ denotes the outward normal.  

%\subsection{Ionic model}\label{sec:ionic_model}

The ionic current $I_\text{ion}$ models the total transmembrane currents through a set of gating and auxiliary variables, respectively denoted by ${\bm z}_E$ and ${\bm z}_S$. The gating variables ${\bm z}_E$ are probabilities that control the opening and closing of ion channels in the cell membrane, whereas the remaining variables ${\bm z}_S$ represent other quantities, such as ion concentration, or intracellular calcium dynamic.  Many different ionic models exist in the literature, depending on the type of cell and the pathophysiological state, but in most cases the overall structure is fixed as in \cref{eq:monodomain} and follows the so-called Hodgkin--Huxley formalism. Specifically, the total current $I_\text{ion}$ is the sum of several ionic currents, each depending on its own set of gating and auxiliary variables. An archetypal instance is the celebrated Hodgkin--Huxley model for electric propagation in axons:
$$
I_\text{ion}(V_m,{\bm z}_E) = g_\text{Na} m^3 h (V_m - V_\text{Na}) + g_\text{K} n^4 (V_m - V_\text{K}) + g_\text{leak} (V_m - V_\text{leak}),
$$
where $\bm{z}_E = ( m, h, n )$ and there is no auxiliary variable ${\bm z}_S$. State-of-the-art cardiac ionic models are more advanced in the sense that they include several more currents and the intracellular calcium dynamic, thus significantly increasing the computational complexity and temporal stiffness. In this respect, it is common to take advantage of the Hodgkin--Huxley structure of the system to avoid implicit numerical schemes or excessively small time steps. In fact, \cref{eq:monodomain_gate} can be rewritten as:
\begin{equation}
\frac{\partial \bm{z}_E}{\partial t} = \bm{\Lambda}_E(V_m)\bigl(\bm{z}_E - {\bm z}_{\infty}(V_m) \bigr), 
%g_E(V_m,z_E)=\Lambda_E(V_m)(z_E-z_{E,\infty}(V_m)),
\end{equation}
where $\bm{\Lambda}_E(V_m) = -(\bm{\alpha}(V_m) + \bm{\beta}(V_m))$ is a diagonal matrix and ${\bm z}_{E,\infty}(V_m) = -\bm{\alpha}(V_m)(\bm{\alpha}(V_m) + \bm{\beta}(V_m))^{-1}$ is the steady state value of ${\bm z}_E$. In general, $\bm{\Lambda}_E(V_m)$ is very stiff for a physiological range of $V_m$, but it can be efficiently integrated with exponential methods (hence the ``E'' notation) since $\bm{\Lambda}_E(V_m)$ is diagonal. In contrast, $\bm{g}_S(V_m,\bm{z}_E,\bm{z}_S)$ is not stiff (``S'' for ``slow'') and rather expensive to evaluate; thus, it is usually integrated with the explicit Euler scheme. This split approach of integrating ${\bm z}_E$ and ${\bm z}_S$ with an exponential and standard explicit method, respectively, was first proposed in the seminal paper of Rush and Larsen \cite{RuL78} and remains very popular even to this day.

\begin{table}[tb]
\centering
\footnotesize
\begin{tabular}{rlll}
\toprule
\textbf{Variable} & \textbf{Description} & \textbf{Value} & \textbf{Units} \\
\midrule
$V_m(t,\bm x)$ & Transmembrane potential & --- & \si{\milli\volt} \\
$\bm{z}_E(t,\bm x)$ & Vector of gating variables & --- & Probability \\
$\bm{z}_S(t,\bm x)$ & Vector of auxiliary ionic variables & --- & Depends \\
$I_\text{ion}(V_m,\bm{z}_E,\bm{z}_S)$ & Transmembrane current density & --- & \si{\uA\per\square\mm} \\
$I_\text{stim}(t,\bm x)$ & Stimulus current density & --- & \si{\uA\per\square\mm} \\
$\chi$ & surface area-to-volume ratio of cardiac cells & $140$ & \si{\per\mm} \\
$C_m$ & cell membrane electrical capacitance & $0.01$ & \si{\micro\farad\per\square\mm} \\
%$\mathbf{D}(\bm x)$ & electric conductivity tensor & --- & \si{\milli\siemens\per\mm} \\
$\sigma_\ell^i$, $\sigma_t^i$, $\sigma_n^i$ & intracellular electric conductivities & $[0.17,0.019,0.019]$ & \si{\milli\siemens\per\mm} \\
$\sigma_\ell^e$, $\sigma_t^e$, $\sigma_n^e$ & extracellular electric conductivities & $[0.62,0.24,0.24]$ & \si{\milli\siemens\per\mm} \\
$\bm{\alpha}(V_m)$, $\bm{\beta}(V_m)$ & transition rates (diagonal matrices) & --- & \si{\per\ms} \\
\bottomrule
\end{tabular}
\caption{Description of variables and parameters for the monodomain model. The values of the coefficients are from \cite{NKB11}.}
\label{tab:monopar}
\end{table}

\subsection{Space discretization}\label{sec:monodomain_space_disc}
The spatial discretization of the monodomain equations \cref{eq:monodomain}
using conforming piecewise linear finite elements leads to the semi-discrete system
\begin{subequations}\label{eq:sd_sys}
	\begin{align}\label{eq:sd_sys_a}
		C_m\mathbf{M}\frac{\dif \bm V_m}{\dif t} &= \chi^{-1}\mathbf{A}\bm{V}_m +\mathbf{M}\bm{I}_\text{stim}(t) - \mathbf{M}I_\text{ion}(t,\bm{V}_m,\bm{z}_E,\bm{z}_S),\\ \label{eq:sd_sys_b}
		\frac{\dif \bm{z}_E}{\dif t} &= g_E(\bm{V}_m,\bm{z}_E),\\ \label{eq:sd_sys_c}
		\frac{\dif \bm{z}_S}{\dif t} &= g_S(\bm{V}_m,\bm{z}_E,\bm{z}_S),
	\end{align}
\end{subequations}
with appropriate initial conditions and where $\mathbf{A}$ and $\mathbf{M}$ are respectively the stiffness and mass matrices. (For the sake of simplicity, we keep using the vector notation for ${\bm z}_S$ and ${\bm z}_E$.) \Cref{eq:sd_sys_b,eq:sd_sys_c} are evaluated point-wise, due to the absence of differential operators in space in the equations. The nonlinear term $I_\text{ion}$ is evaluated component-wise and then interpolated, whereas $\bm{I}_\text{stim}(t)$ is obtained by linear interpolation of $I_\text{stim}(t,\bm{x})$. This is commonly done in the cardiac modeling community and called Ionic Current Interpolation (ICI), see~\cite{PatMirSouWhi11,PHS16}. 
%Since we first evaluate $I_{ion}$ component-wise and then apply the mass matrix $\bm{M}$, in the terminology of Pathmanathan et al. \cite{PatMirSouWhi11}, we perform ionic current interpolation (ICI).
Alternatively, one could first evaluate $\bm{V}_m$, $\bm{z}_E$, $\bm{z}_S$ at the quadrature nodes and then evaluate $I_\text{ion}$ there, an approach known as State Variable Interpolation (SVI). Here, we opt for ICI due to its lower computational cost, as SVI requires interpolating every state variable at the quadrature nodes, which is expensive for ionic models with a large number $N$ of state variables. Finally, the mass matrix $\mathbf{M}$ is always lumped to take advantage of explicit time integrators.  

The semi-discrete equation \cref{eq:sd_sys} can be written as \cref{eq:FSEeq}
for $y=(\bm{V}_m,\bm{z}_E,\bm{z}_S)^\top$ by defining
%\begin{equation}\label{eq:mono_in_FSE}
%	\begin{aligned}	
%	f_F(t,y)&=((\chi C_m \mathbf{M})^{-1}\mathbf{A}\bm{V}_m,\bm{0},\bm{0})^\top,\\
%	f_S(t,y)&= C_m^{-1}(\bm{I}_{stim}(t)-I_{ion}(t,\bm{V}_m,\bm{z}_E,\bm{z}_S)), \bm{0},g_S(\bm{V}_s,\bm{z}_E,\bm{z}_S))^\top ,\\
%	f_E(t,y)&=(\bm{0},g_E(\bm{V}_m,\bm{z}_E),\bm{0})^\top.
%\end{aligned}
%\end{equation}
\begin{subequations}
\begin{gather}\label{eq:mono_in_FSE}
	%\begin{aligned}
		f_F(t,y)=\begin{pmatrix}
			(\chi C_m \mathbf{M})^{-1}\mathbf{A}\bm{V}_m \\ \bm{0} \\ \bm{0}
		\end{pmatrix}, \qquad
		f_S(t,y)=\begin{pmatrix}
		C_m^{-1}(\bm{I}_{stim}(t)-I_{ion}(t,\bm{V}_m,\bm{z}_E,\bm{z}_S)) \\ \bm{0} \\ g_S(\bm{V}_m,\bm{z}_E,\bm{z}_S)\end{pmatrix}, \\
		f_E(t,y)=\begin{pmatrix}
			\bm{0} \\ g_E(\bm{V}_m,\bm{z}_E) \\ \bm{0}
		\end{pmatrix}
		= \Lambda(y)(y-y_\infty(y)), 
		%\end{aligned}
\end{gather}
with
\begin{equation}
		\Lambda(y)=\begin{pmatrix}
	0 & 0 & 0 \\ 0 & \Lambda_E(\bm{V}_m) & 0 \\ 0 & 0 & 0
\end{pmatrix}, \quad
y_\infty(y) = \begin{pmatrix}
	\bm{0} \\ \bm{z}_{E,\infty}(\bm{V}_m) \\ \bm{0}
\end{pmatrix}.
\end{equation}
\end{subequations}
Note that $f_F$ is stiff but very cheap to evaluate, since it only involves a matrix-vector multiplication thanks to mass-lumping. The $f_S$ term, which contains the most complicated terms of the ionic model, is expensive yet non-stiff since the stiff terms of the ionic model are contained in $g_E$. Clearly, $f_E$ is very stiff but 
straightforward to integrate exponentially, since $\Lambda(y)$ is diagonal.

The time stepping methods presented in the sequel are not limited to the present spatial discretization. Still, the resulting mass matrix must be diagonal, or very cheap to invert, regardless of the spatial discretization used, to preserve the efficiency of any explicit method. For standard conforming triangular or tetrahedral 
finite elements, mass-lumping is available up to a certain polynomial order \cite{COHEN2001, GEEVERS2018}.
Alternative spatial discretizations that lead to diagonal or block diagonal mass matrices at any order are, for instance, higher order finite elements with orthogonal basis functions, spectral finite element methods \cite{AFRICA2023111984}, finite differences, and discontinuous Galerkin methods \cite{OGIERMANN2024116806}.

\ifstandalone
\bibliographystyle{abbrv}
\bibliography{../library}
\fi

\section{Multirate Explicit Stabilized Methods}\label{sec:mES}
In this section, we recall the multirate explicit stabilized mRKC method from \cite{AGR22} for 
\begin{equation}\label{eq:FSeq}
y'=f(y)=f_F(y)+f_S(y),  \qquad y(0)=y_0,
\end{equation}
where $f_F$ is fast (stiff) but cheap, while $f_S$ slow (non-stiff or mildly stiff) yet expensive to evaluate.

To simplify the notation, we restrict the discussion to autonomous problems. Still, the methods proposed here also apply to non-autonomous problems while all listed pseudo-codes apply to the general non-autonomous case.
For the sake of brevity, we concentrate here on algorithmic aspects of the mRKC method from a practical point of view. First, in \cref{sec:ES}, we recall standard explicit stabilized Runge--Kutta--Chebyshev (RKC) methods without any multirate time-stepping strategy. Then, in \cref{sec:mod_eq_avg_f}, we summarize the mathematical derivation and main properties of the mRKC method. The full mRKC Algorithm is then listed 
in \cref{sec:algo_mES}. Finally, in \cref{sec:monodomain_mRKC} we apply the mRKC method to the monodomain model \eqref{eq:monodomain}. 

\subsection{Explicit stabilized methods}\label{sec:ES}
The coefficients of standard explicit Runge--Kutta (RK) methods for
\begin{equation}\label{eq:ode}
y'=f(y),\qquad y(0)=y_0,
\end{equation}
are generally optimized for high accuracy. Those high-order methods, however, suffer from severe stability constraints on the time-step.
In contrast, explicit stabilized (ES) methods fix the order $p$ (usually $p=1,2$) while utilizing an increasing number of stages $s\geq p$ and the corresponding free coefficients to optimize for stability. The resulting stability polynomial, $R_s(z)$, of ES methods then satisfies $|R_s(z)|\leq 1$ for $z\in [-\beta s^2,0]$, where $s$ is the number of stages and $\beta>0$ depends on the method. Therefore, the stability domain grows quadratically with $s$, and hence with the number of function evaluations, so that ES methods remain efficient without sacrificing explicitness. 
Moreover, the number of stages $s$ can easily be adapted every
time-step to accommodate the problem's instantaneous stiffness, measured in terms of the spectral radius $\rho$ of the right-hand side's Jacobian.
%which can be estimated cheaply.
Since ES methods are defined via a recursive relation, they can be implemented using only three vectors, regardless of $s$, and are therefore very efficient memory-wise. 

Various explicit stabilized methods, such as RKC \cite{SSV98}, RKL \cite{Meyer2014}, RKU \cite{mRKCcode} and ROCK \cite{AbM01,Abd02}, are available. Both first- 
and second-order versions exist for
RKC, RKL, and RKU methods, while ROCK methods are second- or fourth-order accurate. 
%As mentioned, the stability domain grows as $\beta s^2$ for all those methods, but for different constants $\beta$. For first-order RKC we have $\beta\approx 2$, while for ROCK2 $\beta=0.82$ and for ROCK4 $\beta=0.35$; hence, $\beta$ decreases as the order increases. 
Here, we shall only consider the classical first-order RKC method from \cref{algo:rkc}, see \cite{Abdulle2015d,HoS80,SSV98,Ver80,Ver96,VHS90} for further details.

\IfStandalone{}{\begin{algorithm}}
\begin{algorithmic}[1]
	\Function{RKC\_Step}{$t_n$, $y_n$, $\Delta t$, $f$, $\rho$}
	\State $\varepsilon=0.05$
	\State $s, \ell_s=\ $\Call{get\_stages}{$\Delta t$, $\rho$, $\varepsilon$}
	\State $y_{n+1} =\ $\Call{RKC\_Iteration}{$t_n$, $y_n$, $\Delta t$, $f$, $s$, $\varepsilon$}
	\State \Return $y_{n+1}$
	\EndFunction
	\Function{get\_stages}{$\Delta t$, $\rho$, $\varepsilon$}
	\State $\beta = 2-4\varepsilon /3$
	\State $s=\lceil\sqrt{\Delta t \,\rho/\beta}\rceil$, $\ell_s=\beta s^2$ \label{algo:rkc_s}
	\State \Return $s$, $\ell_s$
	\EndFunction
	\Function{RKC\_Iteration}{$t$, $y$, $\Delta t$, $f$, $s$, $\varepsilon$}
	\State $\mu_j,\nu_j,\kappa_j,c_j=\ $\Call{get\_coeffs}{$s$, $\varepsilon$}
	\State $g_0 = y$
	\State $g_1=g_0+\mu_1\Delta t f(t,g_0)$
	\For{$j=2,\ldots,s$}
	\State $g_j = \nu_j g_{j-1}+\kappa_j g_{j-2}+\mu_j\Delta t f(t+c_{j-1}\Delta t,g_{j-1})$ \label{algo:rkc_iter}
	\EndFor
	\State \Return $g_s$
	\EndFunction
	\Function{get\_coeffs}{$s$, $\varepsilon$}
	\State $\omega_0=1+\varepsilon/s^2$, $\omega_1=T_s(\omega_0)/T_s'(\omega_0)$
	\State $b_j=1/T_j(\omega_0)$, $j=0,\ldots,s$
	\State $\mu_1=\omega_1/\omega_0$
	\State $c_0=0$, $c_1=\mu_1$
	\For{$j=2,\ldots,s$}
	\State $\mu_j=2\omega_1 b_j/b_{j-1}$, $\nu_j=2\omega_0 b_j/b_{j-1}$,  $\kappa_j=-b_j/b_{j-2}$
	\State $c_j=\nu_j c_{j-1}+\kappa_j c_{j-2}+\mu_j$
	\EndFor
	\State \Return $\mu_j,\nu_j,\kappa_j,c_j$
	\EndFunction	
\end{algorithmic}
\IfStandalone{}{
\caption{One step of the RKC method}
\label{algo:rkc}
\end{algorithm}
}

At each time step, the function \textsc{RKC\_Step} of \cref{algo:rkc} computes the approximate solution $y_{n+1}$ at the new time $t_{n+1}$. Its input arguments are the the solution $y_n$ at the current time $t_n$, the step size $\Delta t$, the right-hand side $f$, and an approximation $\rho$ to the spectral radius of $\partial f/\partial y(t_n,y_n)$. The spectral radius $\rho$ is either known or estimated using any simple method, such as Gershgorin's theorem \cite[pg. 89]{HNW08} or a nonlinear power iteration, as described in \cref{algo:powiter} in \cref{app:pow_it} (see also \cite{Lin72,Lin73,Sha91,shampine1980lipschitz,Ver80}).

In \cref{algo:rkc}, we note that the main iteration in \cref{algo:rkc_iter} can be performed using only three distinct vectors for $g_{j-2}$, $g_{j-1}$, and $g_j$. Moreover, at \cref{algo:rkc_s} the number of stages is proportional to the square root of $\rho$, 
instead of $\rho$, as in standard explicit (RK) methods. The damping parameter $\varepsilon>0$ guarantees that the stability polynomial $R_s(z)$ satisfies $|R_s(z)|\leq 1-\varepsilon$ along the negative real axis for $z\in [-\beta s^2,-\delta]$ and small $\delta>0$. It also slightly extends the stability domain into the imaginary direction thereby ensuring a stable narrow strip along the negative real axis; typically, we set $\varepsilon=0.05$. 
The method's coefficients $\mu,\nu,\kappa,c$ depend on Chebyshev polynomials of the first kind, $T_j(x)$. Indeed, the stability polynomial $R_s(z)$ is a shifted and scaled Chebyshev polynomial.

\subsection{Modified equation and averaged force}\label{sec:mod_eq_avg_f}
Explicit stabilized methods, such as RKC, certainly alleviate the stringent 
stability constraint on the time-step due to the stiff term $f_F$ in \cref{eq:FSeq}. 
Nonetheless, the number $s$ of evaluations of the expensive term $f_S$ still depends on the stiffness of $f_F$, as $\rho$, and hence $s$, are still dictated by the stiffness of $f_F$. To ensure that the number of evaluations of $f_S$ only depends on its own stiffness, and no longer on that of $f_F$, the evaluation of $f_F$ and $f_S$ must be decoupled. 

Following \cite{AGR22}, we therefore introduce the modified equation
\begin{equation}\label{eq:modeq}
y_\eta'=f_\eta(y_\eta),\qquad y_\eta(0)=y_0,
\end{equation}
with a free parameter $\eta>0$. Here, the averaged force $f_\eta:\mathbb{R}^n\rightarrow\mathbb{R}^n$ is defined as
\begin{equation}\label{eq:defavgf}
f_\eta(y)=\frac{1}{\eta}(u(\eta)-y),
\end{equation}
where the auxiliary solution $u(s)$ is defined by the auxiliary equation
\begin{equation}\label{eq:auxp}
u'=f_F(u)+f_S(y) \quad s\in [0,\eta],\qquad u(0)=y.
\end{equation}

We now briefly summarize the main properties of the above modified equation and averaged force; for further details, we refer to \cite{AGR22}.

\paragraph{\textbf{Accuracy.}} The modified equation \cref{eq:modeq} is a first-order approximation to \cref{eq:FSeq} with error $\mathcal{O}(\eta)$, since it follows
from \cref{eq:defavgf,eq:auxp} that
\begin{equation}\label{eq:feta_avg}
    f_\eta(y)=f_S(y)+\frac{1}{\eta}\int_0^\eta f_F(u(s))\dif s=f_S(y)+f_F(y)+\mathcal{O}(\eta).
\end{equation}
In fact, \cref{eq:feta_avg} also implies that $f_\eta$ is a particular average of $f_F+f_S$ along the auxiliary solution $u$. 

\paragraph{\textbf{Stability.}} If the multirate problem \cref{eq:FSeq} is contractive \cite[Chapter IV.12]{HaW02}, then, under some assumptions, \cref{eq:modeq} is also contractive and
\begin{equation}\label{eq:stab_mod_eq}
    \Vert y(t)-y_\eta(t)\Vert\leq C(\eta)\int_0^t e^{\mu (t-s)}\Vert f_F(y(s))+f_S(y(s))\Vert \dif s
\end{equation}
holds, where $y(t)$ is the exact solution to \cref{eq:FSeq}, $\mu<0$, and $C(\eta)\leq 1$ with $C(\eta)=\mathcal{O}(\eta)$ for $\eta\to 0$. Hence, the solutions of \cref{eq:FSeq,eq:modeq} remain close even during long-time integration.

\paragraph{\textbf{Stiffness.}}
The spectral radius $\rho_\eta$ of the Jacobian of $f_\eta$ decreases with increasing $\eta$. Hence, for $\eta$ sufficiently large, it holds 
\begin{equation}\label{eq:stiffness_condition_eta}
    \rho_\eta\leq\rho_S,
\end{equation}
where $\rho_S$ if the spectral radius of the Jacobian of $f_S$. As the stiffness of \cref{eq:modeq} then depends only on the slow term $f_S$, the numerical solution of \cref{eq:modeq} with any explicit method will require fewer right-hand side evaluations than that of \cref{eq:FSeq}.

%\GR{We could remove the next paragraph, since the cost is discussed in the next section more in detail.}
%\paragraph{\textbf{Computational cost.}} In practice, $f_\eta$ is evaluated numerically by approximating \cref{eq:auxp}; however, solving \cref{eq:auxp} is relatively cheap since in general $\eta\ll \Delta t$ and the expensive term $f_S$ is frozen. More details are given in \cref{sec:algo_mES}.

Clearly, there is a trade-off between choosing $\eta$ sufficiently small to preserve accuracy in \cref{eq:feta_avg,eq:stab_mod_eq} yet sufficiently large to satisfy the stiffness condition \cref{eq:stiffness_condition_eta}. Fortunately, both conditions are typically satisfied simultaneously, as \cref{eq:stiffness_condition_eta} is already satisfied for quite small $\eta\leq \Delta t$ (often $\eta\ll \Delta t$), where $\Delta t$ is the step size used to solve the modified equation \cref{eq:modeq}. Then, \cref{eq:modeq} not only is a good approximation to \cref{eq:FSeq} but also much less stiff thanks to \cref{eq:stiffness_condition_eta}.
Discretization of the modified equation \cref{eq:modeq} finally leads to the multirate RKC (mRKC) method, first proposed in \cite{AGR22}.

\subsection{Multirate explicit stabilized method algorithm}\label{sec:algo_mES}
Instead of solving \cref{eq:FSeq} directly, the multirate explicit stabilized (mRKC) method \cite{AGR22} solves the modified problem \cref{eq:modeq} using an RKC method. In doing so, at every evaluation of $f_\eta$, $u(\eta)$ is approximated by solving \cref{eq:auxp} with yet another RKC method. Special care is taken when choosing the number of stages $s$ and $m$ needed for \cref{eq:modeq,eq:auxp}, respectively, and the value of $\eta$. For an accuracy and stability analysis of the method, we refer to \cite{AGR22,Ros20}.

In \cref{algo:mRKC}, we list the mRKC method for the numerical solution of \cref{eq:FSeq}.
The input parameters of \textsc{mRKC\_Step} are $t_n$, $y_n$ and $\Delta t$ (as for \textsc{RKC\_Step}), the two right-hand sides $f_F$, $f_S$ and the spectral radii $\rho_F$, $\rho_S$ of their respective Jacobians.
\IfStandalone{}{\begin{algorithm}}
	\begin{algorithmic}[1]
		\Function{mRKC\_Step}{$t_n$, $y_n$, $\Delta t$, $f_F$, $f_S$, $\rho_F$, $\rho_S$}
		\State $s,\ \ell_s=\ $\Call{get\_stages}{$\Delta t$, $\rho_S$, $0.05$}
		\State $\eta = 2\Delta t/\ell_s$ \label{algo:mRKC_def_s}
		\Comment{For RKC, $\ell_s=\beta s^2\approx 2s^2$}
		\State $m, \ell_m = \ $\Call{get\_stages}{$\eta$, $\rho_F$, $0.05$}
		\State $\overline f_\eta(t,y) = \ $\Call{Averaged\_Force}{$t$, $y$, $\eta$, $f_F$, $f_S$, $m$}
		\Comment{$\overline f_\eta$ is the approximated averaged force \cref{eq:defavgf}}
		\State $y_{n+1} =\ $\Call{RKC\_Iteration}{$t_n$, $y_n$, $\Delta t$, $\overline f_\eta$, $s$, $0.05$} \label{algo:mRKC_outer}
            \Comment{This call approximates \cref{eq:modeq}}
		\State \Return $y_{n+1}$
		\EndFunction
		\Function{Averaged\_Force}{$t$, $y$, $\eta$, $f_F$, $f_S$, $m$}
		\State $f_u(r,u) = f_F(r,u)+f_S(t,y)$
		\Comment{$f_u$ is the right-hand side of \cref{eq:auxp}, $f_S$ is frozen}
		\State $u_\eta =\ $\Call{RKC\_Iteration}{$t$, $y$, $\eta$, $f_u$, $m$, $0.05$} \label{algo:mRKC_inner}
        \Comment{This call approximates \cref{eq:auxp}}
		\State \Return $\tfrac{1}{\eta}(u_\eta-y)$
		\EndFunction
	\end{algorithmic}
	\IfStandalone{}{
		\caption{One step of the mRKC method}
		\label{algo:mRKC}
	\end{algorithm}
}
A single step of \cref{algo:mRKC} mainly consists of an outer \textsc{RKC\_Iteration} (\cref{algo:mRKC_outer}), which solves \cref{eq:modeq} with the approximate averaged force $\overline f_\eta$, and an inner \textsc{RKC\_Iteration} (\cref{algo:mRKC_inner}), which solves \cref{eq:auxp} whenever $\overline f_\eta$ is needed. 

Since the mRKC method is composed of two embedded RKC iterations, the stabilization process is ``multiplicative'', in the sense that the number of outer RKC iterations or stages, $s$, suffices to stabilize the $f_S$ term, whereas the product of outer and inner RKC iterations, $s\times m$, suffices to stabilize the $f_F$ term. This property is reflected in the computational cost.

\paragraph{\textbf{Computational cost.}}
For each evaluation of $\overline f_\eta$, there is a call to the inner loop, which involves one evaluation of $f_S$ and $m$ evaluations of $f_F$. Since there are $s$ evaluations of $\overline f_\eta$ during each time step, there are $s$ evaluations of $f_S$ and $s\times m$ evaluations of $f_F$. According to \cref{algo:mRKC}, we compute
\begin{align}\label{eq:s_approx}
\# f_S \text{ eval.}&=s \approx \sqrt{\frac{\Delta t\rho_S}{\beta}},\\ \label{eq:sm_approx}
\# f_F \text{ eval.}&=s\times m \approx \sqrt{\frac{\Delta t\rho_S}{\beta}}\sqrt{\frac{\eta\rho_F}{\beta}} =  \sqrt{\frac{\Delta t\rho_S}{\beta}}\sqrt{\frac{2\Delta t\rho_F}{\beta^2 s^2}}\approx \sqrt{\frac{2\Delta t\rho_F}{\beta^2}}\approx \sqrt{\frac{\Delta t\rho_F}{\beta}}.
\end{align}
Hence, the number of $f_S$ and $f_F$ evaluations depends only on their own inherent stiffness. 

The total number of evaluations of $f_F$ is approximately the same as with a standard RKC method, for which the number of stages is given by $\sqrt{\Delta t\rho/\beta}$ with $\rho\approx\rho_F$. In contrast, the number of evaluations of $f_S$ are drastically reduced. Finally, the ``multiplicative'' property of the mRKC method is clearly apparent in \cref{eq:sm_approx}, as $s\times m$ satisfies precisely the stability condition $\Delta t\rho_F\approx \beta (sm)^2$ for $f_F$.

%\GR{This remark is not very relevant, we could remove it.}
%\begin{rmrk}
%Let us comment on the differences with respect to a Lie--Trotter splitting method for \cref{eq:FSeq}, since in that case the evaluations of $f_S$ would also depend on its own stiffness only. In fact, the method would first solve $y_1'=f_S(y_1)$ and then $y_2'=f_F(y_2)$, both with RKC methods where the number of stages would depend on $f_S$ and $f_F$, respectively. However, this approach would reduce accuracy due to the clear-cut decoupling. In contrast, the mRKC method continuously updates the contribution of both terms and maintains a stronger coupling, hence higher accuracy. Moreover, as already said in \cref{sec:mod_eq_avg_f}, \cref{eq:modeq} preserves the contractivity properties of the original problem. Indeed, it is verified in \cite{AGR22,AbR22b} that the relative error between mRKC and RKC is negligible.
%\end{rmrk}

\subsection{The mRKC method for the monodomain model}\label{sec:monodomain_mRKC}
To apply mRKC to the semi-discrete monodomain model \cref{eq:sd_sys} from \cref{sec:monodomain_space_disc}, we must first rewrite it as in \cref{eq:FSeq}.
The monodomain model \cref{eq:sd_sys}, however, fits the general form
\cref{eq:FSEeq}, which involves the additional $f_E$ term from the ionic model. In this section, 
we briefly describe a straightforward but somewhat simplistic and at times inefficient approach to deal with $f_E$ and thus apply mRKC directly to \cref{eq:sd_sys}. 
Then, in \cref{sec:exp_mES}, we shall adapt mRKC to differential equations as in \cref{eq:FSEeq} by judiciously exploiting the particular structure of $f_E$, which
will yield a more efficient solver for \cref{eq:sd_sys}.

To directly apply mRKC to \cref{eq:sd_sys}, we first need to absorb the $f_E$ term from the ionic model into $f_F$ and $f_S$. With $y, f_F, f_S, f_E$ as in \cref{eq:mono_in_FSE}, we define
\begin{subequations}\label{eq:mono_in_FS}
\begin{equation}
\tilde y' =\tilde f_F(t,\tilde y)+\tilde f_S(t,\tilde y),
\end{equation}
with
\begin{equation}\label{eq:mRCK_fFfS_monodomain}
	\tilde f_F(t,y)=f_F(t,y)+f_E^1(t,y),\qquad \tilde f_S(t,y)=f_S(t,y)+f_E^2(t,y).
\end{equation}
\end{subequations}
In doing so, we have split the ionic model as $f_E^1(t,y)+f_E^2(t,y)=f_E(t,y)$, where $f_E^1$ contains only the stiffest components while $f_E^2$ contains all the remaining ones. Now, mRKC can be applied to the monodomain model \cref{eq:sd_sys} 
by solving \cref{eq:mono_in_FS}.

Depending on the ionic model and the mesh size, the term $\tilde f_F$ in \cref{eq:mono_in_FS} is usually stiffer than $f_F$, but still cheap to evaluate since $f_E^1$ involves only a few components.
In contrast, as $f_E^2$ is significantly less stiff than $f_E^1$, the term $\tilde f_S$ is only mildly stiff but more expensive to evaluate. 
Nevertheless, this approach not only remains unsatisfactory, as it does not exploit the special structure of $f_E$, but also cumbersome to implement because the stiffest gating variables in $f_E$ must be identified. In the next section, we shall introduce a different more efficient mRKC approach, which exploits the particular structure of $f_E$ without the need for splitting the gating variables.

\begin{rmrk} \label{rem:splitting}
For the cardiac models considered here, it turns out that the stiffest components $f_E^1$ from the ionic model always reduced to a single variable for the entire simulation. 
To identify that stiffest gating variable, we temporarily set $\tilde f_F(t,y)=f_E^1(t,y)$ thereby ignoring diffusion. Next, selecting one gating variable at a time in $f_E^1$, we run the simulation while evaluating the stiffness of $\tilde f_F$ and $\tilde f_S$ at each time step with the nonlinear power iteration -- see \cref{algo:powiter} in the Appendix. The stiffest gating variable is that with the highest spectral radius for $\tilde f_F$ over the entire simulation.
\end{rmrk}

\ifstandalone
\bibliographystyle{abbrv}
\bibliography{../library}
\fi

\section{Exponential Multirate Explicit Stabilized Methods}\label{sec:exp_mES}

As shown in \cref{sec:monodomain_mRKC}, the mRKC method can be directly applied to \cref{eq:FSEeq} by splitting $f_E$ and absorbing the two terms into $f_F$ and $f_S$. Still, whenever $f_E$ is stiffer than $f_F$, it will induce a higher number of stages $m$ in \cref{algo:mRKC}. To improve on the efficiency, we shall now exploit the special structure \cref{eq:def_fE} of $f_E$, as in the Rush--Larsen method \cite{RuL78}, by using an exponential Euler step to integrate $f_E$.
The combined mRKC / exponential Euler method then yields 
the exponential mRKC (emRKC) method for \cref{eq:FSEeq}, 
and hence for the monodomain equation \cref{eq:sd_sys}.

\subsection{Exponential Euler method}\label{sec:exp_E}
A very efficient method for the numerical solution of
\begin{equation}\label{eq:Eeq}
	y'=f_E(t,y),\qquad y(0)=y_0,
\end{equation}
with $f_E(y)=\Lambda(y)(y-y_\infty(y))$ as in \cref{eq:def_fE} and $\Lambda(y)$ a diagonal matrix, is the exponential Euler method \cite{HoO10}:
\begin{equation}\label{eq:exp_E}
	y_{n+1} = y_n + \Delta t\,\varphi(\Delta t\,\Lambda(y_n))f_E(y_n),
\end{equation}
where
\begin{equation}
	\varphi(z)=\frac{e^z-1}{z}
\end{equation}
is an exponential-like function. In contrast to the explicit Euler method, $y_{n+1}=y_n+\Delta t\,f_E(y_n)$, the exponential Euler method employs the matrix exponential $\varphi(\Delta t\,\Lambda(y_n))$, which dampens all large eigenvalues (in magnitude) and stabilizes the discrete dynamical system. Hence, the exponential Euler method is not only unconditionally stable but here also relatively cheap because $\Lambda(y)$ is diagonal. Higher-order variants of this method exist, we refer to \cite{HoO10} for a review article.

\subsection{Exponential Euler step for the auxiliary equation}\label{sec:exp_step}
Instead of splitting $f_E$, as in \cref{sec:monodomain_mRKC}, we could also 
entirely include it into the ``fast'' term $f_F$.
%If we would apply the methodology of \cref{sec:mod_eq_avg_f} to \cref{eq:FSEeq}, the %$f_E$ term, due to its stiffness, would be considered to be part of  
Hence, the auxiliary equation \cref{eq:auxp} would read
\begin{equation}\label{eq:auxp_fE}
	u'=f_F(u)+f_S(y)+f_E(u) \quad t\in [0,\eta],\qquad u(0)=y.
\end{equation}
As in \cref{sec:monodomain_mRKC}, however, the cost would again increase whenever $f_E$ is stiffer than $f_F$.

%if then we would solve this problem with an explicit method, it might become expensive %due to the stiffness of $f_E$, which can be higher than $f_F$. 

To design an explicit multirate method for \cref{eq:FSEeq}, which avoids any further stability restriction due to $f_E$ while taking advantage of its particular diagonal structure, we now modify the auxiliary problem \cref{eq:auxp_fE} by taking inspiration from splitting methods as follows: first, advance the $f_E$ term with the exponential Euler method and then integrate the remaining terms. This yields the split auxiliary equation:
\begin{subequations}\label{eq:auxp_exp}
	\begin{align}\label{eq:auxp_exp_a}
		y_E &= y+\eta\,\varphi(\eta\Lambda(y))f_E(y),\\ \label{eq:auxp_exp_b}
		u' &=f_F(u)+f_S(y_E) \quad t\in [0,\eta],\qquad u(0)=y_E.
	\end{align}
\end{subequations}
In contrast to \cref{eq:auxp_fE}, the stiffness of \cref{eq:auxp_exp_b}
is independent from the $f_E$ term; hence, approximating $u$ with any explicit method is cheaper.

The emRKC method is derived similarly to mRKC though by considering the split auxiliary equation \cref{eq:auxp_exp} instead of \cref{eq:auxp_fE}. Hence, emRKC also solves the modified equation \cref{eq:modeq} instead of the original problem \cref{eq:FSEeq} and the definition \cref{eq:defavgf} of the averaged force $f_\eta$ remains identical, except that $u(\eta)$ in \cref{eq:defavgf} is now computed from \cref{eq:auxp_exp}.

\begin{rmrk}\label{rem:auxp_exp_progr}
Alternatively, we may also split the auxiliary equation with $y_E$ as in \cref{eq:auxp_exp_a} yet with \cref{eq:auxp_exp_b} replaced by $u'=f_F(u)+f_S(y_E)+\varphi(\eta\Lambda(y))f_E(y)$ with $u(0)=y$; hence, the contribution of $\varphi(\eta\Lambda(y))f_E(y)$ is added progressively during the integration. In our numerical experiments, both splittings led to a similar performance.
\end{rmrk}

\subsection{Pseudo-code for the exponential multirate explicit stabilized method}\label{sec:algo_exp_mES}
In \cref{algo:emRKC}, we list a pseudo-code for the emRKC method. Its main differences with respect to \cref{algo:mRKC} are: first, the additional input argument, $f_E$, for \textsc{emRKC\_Step} and \textsc{Averaged\_Force}; second, the redefinition of the $f_u$ term at \cref{algo:emRKC_def_yE,algo:emRKC_redef_fu}; third, the different input argument in \cref{algo:emRKC_inner}. To avoid round-off errors due to small diagonal values in $\Lambda(y)$ or small $\eta$, the computation of the exponential Euler step is done employing the identity indicated in \cref{algo:emRKC_expE_id}.

\IfStandalone{}{\begin{algorithm}}
	\begin{algorithmic}[1]
		\Function{emRKC\_Step}{$t_n$, $y_n$, $\Delta t$, $f_F$, $f_S$, $f_E$, $\rho_F$, $\rho_S$}
		\State $s,\ \ell_s=\ $\Call{get\_stages}{$\Delta t$, $\rho_S$, $0.05$}
		\State $\eta = 2\Delta t/\ell_s$ \label{algo:emRKC_def_s}
		\State $m, \ell_m = \ $\Call{get\_stages}{$\eta$, $\rho_F$, $0.05$}
		\State $\overline f_\eta(t,y) = \ $\Call{Averaged\_Force}{$t$, $y$, $\eta$, $f_F$, $f_S$, $f_E$, $m$}
		\State $y_{n+1} =\ $\Call{RKC\_Iteration}{$t_n$, $y_n$, $\Delta t$, $\overline f_\eta$, $s$, $0.05$} \label{algo:emRKC_outer}
		\State \Return $y_{n+1}$
		\EndFunction
		\Function{Averaged\_Force}{$t$, $y$, $\eta$, $f_F$, $f_S$, $f_E$, $m$}
		\State $y_E =\ $\Call{Exponential\_Euler}{$y$, $\eta$, $\Lambda(y)$, $y_\infty(y)$} \label{algo:emRKC_def_yE}
		\State $f_u(r,u) = f_F(r,u)+f_S(t,y_E)$ \label{algo:emRKC_redef_fu}
		\State $u_\eta =\ $\Call{RKC\_Iteration}{$t$, $y_E$, $\eta$, $f_u$, $m$, $0.05$} \label{algo:emRKC_inner}
		\State \Return $\tfrac{1}{\eta}(u_\eta-y)$
		\EndFunction
		\Function{Exponential\_Euler}{$y$, $\eta$, $\Lambda$, $y_\infty$}
		\State \Return $y +(e^{\eta\Lambda}-I)(y-y_\infty)$
		\Comment{$\eta\varphi(\eta\Lambda)\Lambda(y-y_\infty) = (e^{\eta\Lambda}-I)(y-y_\infty)$} \label{algo:emRKC_expE_id}
		\EndFunction
	\end{algorithmic}
	\IfStandalone{}{
		\caption{One step of the emRKC method}
		\label{algo:emRKC}
	\end{algorithm}
}

\subsection{Stability analysis}\label{sec:stab_analysis}
Here we analyse the stability of the emRKC method starting from previous results for mRKC. In particular, we show that there are no constraints on the step size $\Delta t$ for sufficiently many internal stages $s,m$.
Hence, we consider the standard linear scalar test equation,
\begin{equation}\label{eq:test_eq_emRKC}
	y'=\lambda_F y+\lambda_S y+\lambda_E y,\qquad y(0)=y_0,
\end{equation}
where $\lambda_F,\lambda_S,\lambda_E\leq 0$ correspond to $f_F,f_S,f_E$, respectively. Before tackling the analysis for the emRKC method, we briefly recall the stability analysis for the mRKC method \cite{AGR22,Ros20}.

\subsubsection{Stability analysis for the mRKC method}
The stability analysis of the mRKC scheme \cite{AGR22,Ros20} also relies on \cref{eq:test_eq_emRKC} yet with $\lambda_E=0$; hence, we consider the test equation
\begin{equation}\label{eq:test_eq_mRKC}
	y'=\lambda_F y+\lambda_S y,\qquad y(0)=y_0.
\end{equation}
First, we recall a key result from \cite{AGR22}, crucial for proving the stability of the mRKC method when applied to the test equation \cref{eq:test_eq_mRKC}.
\begin{lemma}\label{lemma:stab_mRKC}
	Let $P_m(z)$ be the stability polynomial of the inner RKC iteration of mRKC (\cref{algo:mRKC_inner} of \cref{algo:mRKC}). Let $\Phi_m(z)=(P_m(z)-1)/z$. If
	\begin{equation}\label{eq:cond_stab_mRKC_raw}
		-\ell_s\leq \Delta t \,\Phi_m(\eta\lambda_F)(\lambda_F+\lambda_S)\leq 0
	\end{equation}
	holds, the mRKC method applied to the test equation \cref{eq:test_eq_mRKC} is stable.
\end{lemma}
\begin{proof}[Sketch of proof.] 
Since mRKC is equivalent to an $s$-stage RKC method applied to the averaged force $\overline f_\eta(y)$ (\cref{algo:mRKC_outer} of \cref{algo:mRKC}), we conclude that mRKC is stable if
\begin{equation}\label{eq:cond_stab_mRKC_raw_equiv}
	-\ell_s\leq \Delta t \frac{\partial\overline f_\eta}{\partial y}(y)\leq 0.
\end{equation}
%which are the stability conditions for an $s$-stage RKC method applied to $\overline %f_\eta(y)$.
To complete the proof, one shows that $\overline f_\eta(y)=\Phi_m(\eta\lambda_F)(\lambda_F+\lambda_S)y$ when mRKC is applied to the test equation \cref{eq:test_eq_mRKC}; thus, condition \cref{eq:cond_stab_mRKC_raw} is equivalent to \cref{eq:cond_stab_mRKC_raw_equiv}.
\end{proof}
In \cite{Ros20}, the conditions on $s,m,\eta$ to satisfy \cref{eq:cond_stab_mRKC_raw} were studied in two different scenarios.
In the first case, $\lambda_F,\lambda_S$ can take any non-positive values and general stability conditions on $s,m,\eta$ are derived. In the second case, it is assumed that $\lambda_F\ll\lambda_S\leq 0$, thus yielding milder stability conditions. We recall these conditions in \cref{lemma:stab_cond} below.

\begin{lemma}\label{lemma:stab_cond}
	Let $\rho_F=|\lambda_F|$ and $\rho_S=|\lambda_S|$, $\ell_s=\beta s^2$ and $\ell_m=\beta m^2$. Define the hypotheses:
	\begin{myenum}
		\item: $\lambda_F,\lambda_S\leq 0$ and $s,m,\eta$ such that
		\begin{equation}\label{eq:stab_cond_A}
			\Delta t\, \rho_S \leq \ell_s , \qquad \eta\,\rho_F\leq \ell_m, \qquad \eta=\frac{6\Delta t}{\ell_s}\frac{m^2}{m^2-1}.
		\end{equation}\label{item:cond_A}
		\item: $\lambda_F\ll\lambda_S\leq 0$ and $s,m,\eta$ such that
		\begin{equation}\label{eq:stab_cond_B}
			\Delta t\, \rho_S \leq \ell_s , \qquad \eta\,\rho_F\leq \ell_m, \qquad \eta=\frac{2\Delta t}{\ell_s}.
		\end{equation}\label{item:cond_B}		
	\end{myenum}	
	If either \cref{item:cond_A} or \cref{item:cond_B} are true, then \cref{eq:cond_stab_mRKC_raw} is satisfied.
\end{lemma}
For the proof of this lemma and more details on the condition $\lambda_F\ll\lambda_S\leq 0$\footnote{Essentially we need $\lambda_F\leq 4\lambda_S$.} we refer to \cite{Ros20}. By combining \cref{lemma:stab_mRKC,lemma:stab_cond}, we eventually prove the stability of the mRKC method for the test equation.
\begin{theorem}\label{thm:stab_mRKC}
	If either hypothesis \cref{item:cond_A} or \cref{item:cond_B} in \cref{lemma:stab_cond} is satisfied, the mRKC method applied to the test equation \cref{eq:test_eq_mRKC} is stable.
\end{theorem}
\begin{proof}
	Follows from \cref{lemma:stab_mRKC,lemma:stab_cond}.
\end{proof}
Note that \cref{item:cond_B} places stronger restrictions on $\lambda_F,\lambda_S$, but also yields a smaller $\eta$ than \cref{item:cond_A}, and hence a smaller number of stages $m$ which results in fewer $f_F$ function evaluations and therefore a more efficient algorithm.
For the context of this work, it makes sense to assume that \cref{item:cond_B} is satisfied, hence $\lambda_F\ll\lambda_S\leq 0$, which stands for $f_F$ significantly stiffer than $f_S$. Indeed, it was found in \cite{AGR22}, that in general for systems of ODEs stemming from the spatial discretization of parabolic PDEs, the stability conditions \cref{eq:stab_cond_B} yield stable solutions.
Note that \cref{eq:stab_cond_B} are the constraints encoded in \cref{algo:mRKC,algo:emRKC}.

We also remark that there is no constraint on the step size $\Delta t$. Indeed, given $\rho_F,\rho_S$, one can always choose $\Delta t$ according to the desired accuracy. Stability is then achieved by setting $s$ sufficiently large to ensure $\Delta t \rho_S \leq \ell_s$ and finally computing $\eta$ and $m$.

\subsubsection{Stability analysis for the emRKC method}
The stability analysis of the emRKC method follows from that for the mRKC method. In particular, we find that the same stability conditions hold and therefore there are no constraints on the step size $\Delta t$.
\begin{theorem}\label{thm:stab_emRKC}
If either hypothesis \cref{item:cond_A} or \cref{item:cond_B} in \cref{lemma:stab_cond} is satisfied, the emRKC method applied to the test equation \cref{eq:test_eq_emRKC} is stable.
\end{theorem}
\begin{proof}
Following \cite{AGR22}, we deduce for emRKC applied to \cref{eq:test_eq_emRKC} that
\begin{equation}
	u_\eta = (P_m(\eta\lambda_F)+\Phi_m(\eta\lambda_F)\eta\lambda_S)y_E=(P_m(\eta\lambda_F)+\Phi_m(\eta\lambda_F)\eta\lambda_S)e^{\eta\lambda_E}y.
\end{equation}
The averaged force $\overline f_\eta(y)$ is therefore given by
\begin{equation}
\begin{aligned}
\overline f_\eta(y)&=\frac{1}{\eta}(u_\eta-y)= \frac{1}{\eta}\left( (P_m(\eta\lambda_F)+\Phi_m(\eta\lambda_F)\eta\lambda_S)e^{\eta\lambda_E}-1  \right)y\\
&=\frac{1}{\eta}\left( (P_m(\eta\lambda_F)-1+\Phi_m(\eta\lambda_F)\eta\lambda_S)e^{\eta\lambda_E}+e^{\eta\lambda_E}-1  \right)y\\
&=\left( \Phi_m(\eta\lambda_F)(\lambda_F+\lambda_S)e^{\eta\lambda_E}+\varphi(\eta\lambda_E)\lambda_E \right)y.
\end{aligned}
\end{equation}
Since emRKC is an $s$-stage RKC method applied to $\overline f_\eta(y)$ (\cref{algo:emRKC_outer} of \cref{algo:emRKC}), it is stable if
\begin{equation}\label{eq:cond_lambaeta}
	-\ell_s\leq \Delta t  \frac{\partial \overline f_\eta}{\partial y}(y) \leq 0.
\end{equation}
By using $\varphi(z)\geq 0$, $\lambda_E\leq 0$ and \cref{eq:cond_stab_mRKC_raw}, we obtain
\begin{equation}
\begin{aligned}
	0&\geq \Delta t  \frac{\partial \overline f_\eta}{\partial y}(y) = \Delta t\left( \Phi_m(\eta\lambda_F)(\lambda_F+\lambda_S)e^{\eta\lambda_E}+\varphi(\eta\lambda_E)\lambda_E \right)\\
	&\geq -\ell_s e^{\eta\lambda_E}+\frac{\Delta t}{\eta}(e^{\eta\lambda_E}-1)=-\ell_s + (e^{\eta\lambda_E}-1)\left(\frac{\Delta t}{\eta}-\ell_s\right)\geq -\ell_s.
\end{aligned}
\end{equation}
The last inequality follows from $e^{\eta\lambda_E}-1\leq 0$ and the fact that $\frac{\Delta t}{\eta}-\ell_s\leq 0$ holds for both choices of $\eta$ in \cref{eq:stab_cond_A} or \cref{eq:stab_cond_B}.
\end{proof}
Note that $\lambda_E$ does not impact the choice of $s$,$m$, or $\eta$ due to its exponential integration.

\subsection{Accuracy analysis}\label{sec:acc_analysis}
To prove that the emRKC method is first-order accurate, we assume that $f_F,f_S,f_E$ are sufficiently smooth. Furthermore, for simplicity of presentation, we restrict ourselves to the autonomous case when $f_F,f_S,f_E$ are time independent.

\begin{theorem}\label{thm:convergence}
The emRKC method is first-order accurate.
\end{theorem}
\begin{proof}
We start by computing a Taylor expansion of the exact auxiliary solution $u(\eta)$ defined by the split auxiliary equation \cref{eq:auxp_exp}. 
We have $y_E=y+\eta f_E(y)+\mathcal{O}(\eta^2)$ and $u(s)=u(0)+\mathcal{O}(s)$. Since $u(0)=y_E$, we thus obtain
\begin{equation}\label{eq:acc_1}
	\begin{aligned}
		u(\eta)&=u(0)+\int_0^\eta u'(s)\dif s=y_E+\int_0^\eta f_F(u(s))+f_S(y_E)\dif s=y_E+\int_0^\eta f_F(y_E)+f_S(y_E)+\mathcal{O}(s)\dif s\\
		&=y_E+\eta f_F(y_E)+\eta f_S(y_E)+\mathcal{O}(\eta^2)=y+\eta (f_F(y)+f_S(y)+f_E(y))+\mathcal{O}(\eta^2).
	\end{aligned}
\end{equation}
Since the inner RKC method (\cref{algo:emRKC_inner} of \cref{algo:emRKC}) used to approximate the integral in \cref{eq:acc_1} is first-order accurate, it introduces a local error of $\mathcal{O}(\eta^2)$; hence, the numerical auxiliary solution $u_\eta$ satisfies
\begin{equation}
	u_\eta = u(\eta)+\mathcal{O}(\eta^2)=y+\eta (f_F(y)+f_S(y)+f_E(y))+\mathcal{O}(\eta^2).
\end{equation}
Therefore, the numerical averaged force $\overline f_\eta(y)$ computed in \cref{algo:emRKC} satisfies
\begin{equation}\label{eq:acc_tfe}
	\overline f_\eta(y)=\frac{1}{\eta}(u_\eta-y)=f_F(y)+f_S(y)+f_E(y)+\mathcal{O}(\eta).
\end{equation}
Hence it is a first-order approximation in $\eta$ of the exact right-hand side. At any time step given $y_n$, the outer RKC method (\cref{algo:emRKC_outer} of \cref{algo:emRKC}) computes $y_{n+1}$ as
\begin{equation}\label{eq:acc_2}
	\begin{aligned}
		g_0&=y_n,\quad g_1=g_0+\mu_1\Delta t \overline f_\eta(g_0),\\
		g_j&= \nu_j g_{j-1}+\kappa_j g_{j-2}+\mu_j\Delta t \overline f_\eta(g_{j-1}),\quad j=2,\ldots,s,\\
		y_{n+1}&=g_s.
	\end{aligned}
\end{equation}
By setting $d_j=g_j-y_n$ and using $\nu_j+\kappa_j=1$ we can rewrite \cref{eq:acc_2} as
\begin{equation}
	\begin{aligned}
		d_0&=0,\quad d_1=\mu_1\Delta t \overline f_\eta(t_n,y_n)\\
		d_j&= \nu_j d_{j-1}+\kappa_j d_{j-2}+\mu_j\Delta t \overline f_\eta(y_n+d_{j-1}),\quad j=2,\ldots,s,\\
		y_{n+1}&=y_n+d_s.
	\end{aligned}
\end{equation}
From \cite[Lemma 3.1]{CrR22} with $A=0$ and $r_j=\mu_j\Delta t \overline f_\eta(y_n+d_{j-1})$ we obtain
\begin{equation}
	y_{n+1}=y_n+\sum_{j=1}^s\frac{b_s}{b_j}U_{s-j}(\omega_0)r_j=y_n+\Delta t\sum_{j=1}^s\frac{b_s}{b_j}U_{s-j}(\omega_0)\mu_j  \overline f_\eta(y_n+d_{j-1}),
\end{equation}
where the $b_j$ are as in \cref{algo:rkc} and $U_j(x)$ denotes the Chebyshev polynomial of the second kind of degree $j$.
Since $d_j=\mathcal{O}(\Delta t)$ (shown recursively) it follows
\begin{equation}\label{eq:acc_3}
	y_{n+1}=y_n+\Delta t\sum_{j=1}^s\frac{b_s}{b_j}U_{s-j}(\omega_0)\mu_j \left( \overline f_\eta(y_n)+\mathcal{O}(\Delta t)\right) =y_n+\Delta t \overline f_\eta(y_n)+\mathcal{O}(\Delta t^2),
\end{equation}
where the last equality follows from the identity $\sum_{j=1}^s\frac{b_s}{b_j}U_{s-j}(\omega_0)\mu_j=1$ (see \cite{CrR22}). Finally, we use \cref{eq:acc_3} and \cref{eq:acc_tfe} to obtain
\begin{equation}
	y_{n+1}= y_n+\Delta t (f_F(y_n)+f_S(y_n)+f_E(y_n))+\mathcal{O}(\Delta t^2+\Delta t\eta).
\end{equation}
Since $\eta=\mathcal{O}(\Delta t)$, the local error is $\mathcal{O}(\Delta t^2)$ and the method is first-order accurate.
\end{proof}

\ifstandalone
\bibliographystyle{abbrv}
\bibliography{../library}
\fi

\section{Numerical experiments}\label{sec:num}
In this section, we apply the mRKC and emRKC methods to the monodomain model \cref{eq:monodomain} and compare them with the popular implicit-explicit Rush--Larsen (IMEX-RL) method described below. 
%
%we perform numerical experiments on the monodomain model \cref{eq:monodomain} and %compare the efficiencies of the mRKC and emRKC methods against a baseline method %presented below. 
Their performance depends on a number of factors, such as the mesh size and regularity, the nonlinearity, the stiffness, and the number of ionic state variables; moreover, the performance of the IMEX-RL method also depends on the type of 
solver used during each time step. To take into account all these factors, we perform numerical experiments either with structured or unstructured meshes for different mesh sizes, while considering ionic models of varying size and stiffness, both in two and three space dimensions thereby including both direct and iterative linear solvers into the comparison. 

More specifically, we first assess in  \cref{sec:eff_stab} the effective stability properties of the emRKC method when compared with the theoretical results of \cref{sec:stab_analysis}. 
Then, in \cref{sec:cuboid2d} we perform a two-dimensional experiment on a structured mesh, where we compare the various methods for different ionic models and mesh sizes. In \cref{sec:cuboid3d}, we repeat the experiment yet for a three-dimensional structured mesh; hence, an iterative solver is used for the linear systems needed in the IMEX-RL method. Next in \cref{sec:fastl_LA}, we consider a more realistic three-dimensional geometry of the left atrium with an unstructured mesh. 
Again we repeat the numerical experiment 
for the same geometry but in the presence of unhealthy fibrotic tissue when the diffusion tensor is discontinuous.  We conclude in \cref{sec:scalability} with a scalability experiment.

\paragraph{\textbf{The implicit-explicit-Rush--Larsen method}}\label{sec:imexexp}
The well-known IMEX-RL method for the monodomain model will be used as a benchmark to assess the performance of mRKC and emRKC. It is based on an IMEX scheme for \cref{eq:sd_sys_a} and the Rush--Larsen approach \cite{RuL78} for the ionic model \cref{eq:sd_sys_b,eq:sd_sys_c}; hence, $\bm{z}_S$ is integrated explicitly whereas the gating variables $\bm{z}_E$ via an exponential integrator. 
%Henceforth, we refer to this method as IMEX-RL. 

The IMEX-RL method proceeds as follows.
First, it solves the ODE subsystem \cref{eq:sd_sys_b,eq:sd_sys_c} with the Rush--Larsen method \cite{RuL78}. Hence, at the beginning of each time step, the method freezes the electric potential $\bm V_m$ and solves \cref{eq:sd_sys_b} with the exponential Euler method. Then it solves \cref{eq:sd_sys_c} with the explicit Euler method. Finally, it solves \cref{eq:sd_sys_a} by treating the linear term implicitly and the remaining terms explicitly. A pseudo-code for the resulting IMEX-RL method is provided in \cref{algo:imexexp}.

Since the stiffest terms $g_E(\bm{V}_m,\bm{z}_E)$ and $\mathbf{A}\bm{V}_m$ are treated exponentially and implicitly, respectively, the scheme only has a mild stability constraint. It is also quite efficient because $\Lambda_E(\bm{V}_m)$ is diagonal and the implicit step is linear in the unknown. 

\IfStandalone{}{\begin{algorithm}}
	\begin{algorithmic}[1]
		\Function{IMEX-RL\_Step}{$t_n$, $\bm{V}_{m,n}$, $\bm{z}_{E,n}$, $\bm{z}_{S,n}$, $\Delta t$}
		\State $\bm{z}_{E,n+1} = \bm{z}_{E,n}+(e^{\Delta t\Lambda_E(\bm{V}_{m,n})}-I)(\bm{z}_{E,n}-\bm{z}_{E,\infty}(\bm{V}_{m,n}))$ \label{algo:imexexp_exp}
		\State $\bm{z}_{S,n+1} = \bm{z}_{S,n}+\Delta t g_S(\bm{V}_{m,n},\bm{z}_{E,n+1},\bm{z}_{S,n})$ \label{algo:imexexp_slow}
		\State Solve for $\bm{V}_{m,n+1}$:
		$C_m\mathbf{M}\frac{\bm{V}_{m,n+1}-\bm{V}_{m,n}}{\Delta t}=\chi^-1\mathbf{A}\bm{V}_{m,n+1}+\mathbf{M}\bm{I}_{stim}(t_n)-\mathbf{M} I_{ion}(t_n,\bm{V}_{m,n},\bm{z}_{E,n+1},\bm{z}_{S,n+1})$.\label{algo:imexexp_implicit}
		\State \Return $\bm{V}_{m,n+1}$, $\bm{z}_{E,n+1}$, $\bm{z}_{S,n+1}$
		\EndFunction
	\end{algorithmic}
	\IfStandalone{}{
		\caption{One step of the IMEX-RL method.}
		\label{algo:imexexp}
	\end{algorithm}
}

\paragraph{\textbf{Computational setup}} 
The code utilized for conducting the numerical experiments in this section is publicly accessible in the repository \cite{emRKCcode}.
The monodomain model \cref{eq:monodomain} is implemented in Python employing the FEniCSx \cite{AlnaesEtal2014,BasixJoss,ScroggsEtal2022} library, where, under the hood, all expensive computations are performed in C/C++. In particular, all linear algebra operations are performed through the PETSc \cite{petsc-web-page,DalcinPazKlerCosimo2011} library.
As ionic models should not be implemented by hand, due to their high complexity, their definitions were obtained from the CellML \cite{lloyd2008cellml} model repository; in doing so, we adapted the C code generated by the Myokit \cite{CCL16} library and created wrappers to call it from Python. Numerical experiments of \cref{sec:eff_stab,sec:cuboid2d} are performed on our workstation, while those of \cref{sec:cuboid3d,sec:fastl_LA,sec:fastl_LA_fibrosis,sec:scalability} on the Piz Daint supercomputer at the Swiss National Supercomputing Centre (CSCS).

The mRKC method always incorporates the splitting of the gating variables described in \cref{eq:mRCK_fFfS_monodomain}---see also Remark \ref{rem:splitting}. In contrast, the emRKC method employs the same splitting as the IMEX-RL method, which is the usual one used in the cardiac community.

Unless stated otherwise, the initial values for $\bm{V}_m$, $\bm{z}_E$, $\bm{z}_S$ are provided by the ionic model's definition and are uniformly set throughout the computational domain. Here we shall consider the Hodgkin--Huxley (HH) \cite{HodHux52}, the Courtemanche--Ramirez--Nattel (CRN) \cite{CRN98} and the ten Tusscher--Panfilov (TTP) \cite{tenPan06} model. Thus our numerical experiments span a wide range of situations from the small nonstiff HH model, to the larger and mildly stiff CRN model, and even to the large and very stiff TTP model, see \cref{tab:ionicmodel_sizestiff} and \cite{Spiteri2010,SpD12}.
In all numerical experiments, we use the parameter values specified in \cref{tab:monopar}. Finally, 
the vector fields $\bm a_l(\bm x),\bm a_t(\bm x), \bm a_n(\bm x)$ (see \cref{eq:def_D}), which define the fibers' orientation in \cref{sec:fastl_LA}, are generated with lifex-fiber \cite{AfrPieFedDedQua23,piersanti2021modeling}.
The applied stimulus $\bm{I}_\text{stim}(t)$ is specified individually for each experiment. 

\begin{table}[htb]
	\begin{center}
		\begin{tabular}{rrrr}
  \toprule
			 & HH & CRN & TTP \\ \midrule		
    $N$ & 3 & 20 & 18 \\
    $\rho_{\max}$ & 40 & 130 & 1000 \\
    \bottomrule
		\end{tabular}
	\end{center}
	\caption{Ionic models' size and stiffness, with $N$ the number of state variables and $\rho_{\max}$ the maximal spectral radius of the right-hand side's Jacobian in \cref{eq:sd_sys}.}
	\label{tab:ionicmodel_sizestiff}
\end{table}

The linear systems needed in the IMEX-RL method are solved differently in two or three space dimensions: in 2--D we use the Cholesky factorization, and in 3--D the conjugate gradient method preconditioned using the algebraic multigrid library hypre \cite{falgout2002hypre}. For mRKC and emRKC, the spectral radii $\rho_F,\rho_S$ are computed only once at the beginning of the simulation employing \cref{algo:powiter} in \cref{app:pow_it}.

During the numerical experiments, we always monitor the relative error in the potential $\bm{V}$ at the final time with respect to the $L^2(\Omega)$-norm:
\begin{equation}
\text{rel. }L^2\text{-err.}=\frac{\Vert \bm{V^\star}-\bm{V}_N\Vert_{L^2(\Omega)} }{\Vert \bm{V^\star} \Vert_{L^2(\Omega)}},
\end{equation}
where $\bm{V}_N$ is the numerical solution obtained with either one of the IMEX-RL, mRKC, emRKC methods and $\bm{V^\star}$ is a reference solution to \cref{eq:sd_sys} computed with a much smaller time-step $\Delta t=\SI{e-4}{\milli\second}$ using a method similar to IMEX-RL but fully explicit, i.e., \cref{algo:imexexp_implicit} in \cref{algo:imexexp} is performed explicitly; we denote this method EXEX-RL. The final time is chosen such that a propagating wave is located inside the domain and the solution is still far from equilibrium.

\subsection{Practical stability limit}\label{sec:eff_stab}
In \cref{sec:stab_analysis} we proved for a scalar linear test equation that the emRKC method is stable for all step sizes $\Delta t>0$ if the numbers of stages $s,m$ are sufficiently large and satisfy the stability condition \cref{eq:stab_cond_B}. In this first numerical experiment, we study the effective stability properties of emRKC in a more realistic setting. 

To assess the practical stability properties of the emRKC method, we solve the monodomain equation \cref{eq:monodomain} with increasing step size $\Delta t$ and verify that the norm of the solution remains bounded, ignoring accuracy here. To push the method to its stability limit, we choose a setting where the stiffness of $f_F$ and $f_E$ is severe, i.e., we choose a relatively small mesh size $\Delta x=\SI{0.025}{\milli\metre}$ in a one dimensional domain $\Omega=[0,20]\si{\milli\metre}$ and employ the TTP ionic model. Although the stiffness of $f_S$ depends on the ionic model, too, it usually does only mildly so. We integrate the monodomain equation up to the final time $T=\SI{100}{\milli\second}$ with an applied stimulus $\chi I_\text{stim}(t,\bm x)=\SI[per-mode=symbol]{50}{\micro\ampere\per\milli\metre\cubed}$ for $\bm x\in [0,1.5]\si{\milli\metre}$ and $t\leq \SI{2}{\milli\second}$, and else $\chi I_\text{stim}(t,\bm x)=0$.

\begin{figure}
	\begin{subfigure}[t]{\subfigsized\textwidth}
		\centering
		\includegraphics[scale=\plotimscalet,trim={0 3mm 0 0},clip]{\currfiledir 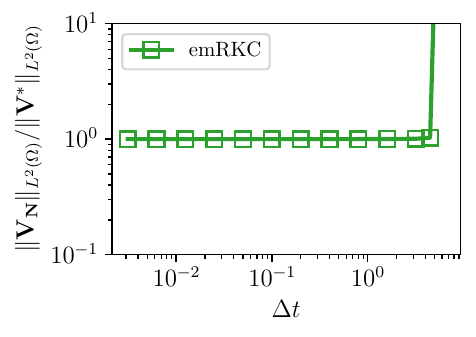}
		\caption{Relative norm of $\mathbf{V}_N$ as $\Delta t$ increases.}
		\label{fig:cub1d_stab_norm}
	\end{subfigure}
	\begin{subfigure}[t]{\subfigsized\textwidth}
		\centering
		\includegraphics[scale=\plotimscalet,trim={0 3mm 0 0},clip]{\currfiledir 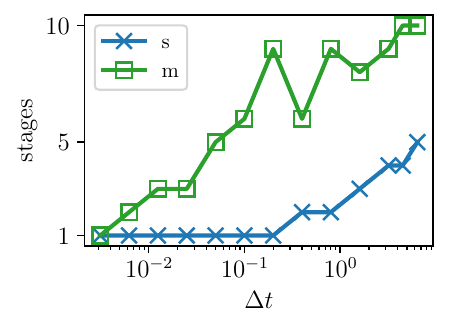}
		\caption{Number of stages $s,m$ as $\Delta t$ increases.}
		\label{fig:cub1d_stab_stages}
	\end{subfigure}
	\caption{Effective stability. Relative norm of the numerical solution $\mathbf{V}_N$ and number of stages for increasing step size $\Delta t$. Recall that $s$, $m$ are the number of stages used by the outer and inner RKC iterations, respectively.}
	\label{fig:cub1d_stab}
\end{figure}

In \cref{fig:cub1d_stab_norm} we display the relative norm of the numerical solution $\mathbf{V}_N$ with respect to a reference solution $\mathbf{V}^*$ as a function of the step size $\Delta t$. We see that the method remains stable up to $\Delta t\leq \SI{4.5}{\milli\second}$. In contrast, a standard explicit scheme, such as the EXEX-RL method described at the beginning of \cref{sec:num} (i.e. as in \cref{algo:imexexp} but performing \cref{algo:imexexp_implicit} explicitly), remains stable only up to $\Delta t\leq \SI{0.003}{\milli\second}$. In \cref{fig:cub1d_stab_stages} we show the number of stages $s,m$ taken by emRKC. We observe that $m>1$ already for small $\Delta t$, indicating that the diffusion $f_F$ has to be stabilized very soon by the inner iterations of emRKC. The outer stages are needed to stabilize the non stiff term $f_S$, as indicated by $s>1$ only for relatively large $\Delta t$. The ``V'' shaped behavior of $m$ is typical and further discussed in the next \cref{sec:cuboid2d}.

In \cref{sec:stab_analysis}, we proved that for the test equation \cref{eq:test_eq_emRKC} the step size $\Delta t$ can be taken arbitrarily large, as the emRKC method remains stable as long as the number of stages $s,m$ are sufficiently large and satisfy the stability condition \cref{eq:stab_cond_B}. In contrast, as shown in the current experiment, due to nonlinearity the maximal step size is bounded in practice. Nevertheless, this practical bound is remarkably large, $\Delta t\leq \SI{4.5}{\milli\second}$ when compared to the stability bound for a standard explicit method, where $\Delta t\leq \SI{0.003}{\milli\second}$, thus yielding a maximal step size increase by a factor of 1500! In that sense, we claim that in practice emRKC has no stability restriction on the step size, since we rarely use $\Delta t>\SI{0.1}{\milli\second}$ because of accuracy requirements.

\subsection{Two-dimensional rectangle}\label{sec:cuboid2d}
First, we consider the two-dimensional version of the benchmark problem proposed in \cite{NKB11}, hence \cref{eq:monodomain} with $\Omega=[0,20]\times [0,7]\,\si{\milli\metre\squared}$ and the final time $T=\SI{25}{\milli\second}$. To initiate the propagation of the action potential, the applied stimulus is defined as $\chi I_\text{stim}(t,\bm x)=\SI[per-mode=symbol]{50}{\micro\ampere\per\milli\metre\cubed}$ if $t\leq \SI{2}{\milli\second}$ and $\Vert\bm x\Vert_\infty\leq \SI{1.5}{\milli\metre}$, else $\chi I_\text{stim}(t,\bm x)=0$. 
%Here, we fix the end time at $T=\SI{25}{\milli\second}$.

\begin{figure}
	\begin{subfigure}[t]{\subfigsizet\textwidth}
		\centering
		\includegraphics[scale=\plotimscalet,trim={0 3mm 0 0},clip]{\currfiledir 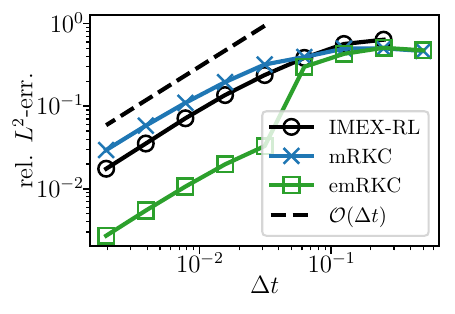}
		\caption{HH model, $\Delta x=0.2 \mm$.}
		\label{fig:cub2d_conv_ref_0_HH}
	\end{subfigure}\hfill%
	\begin{subfigure}[t]{\subfigsizet\textwidth}
		\centering
		\includegraphics[scale=\plotimscalet,trim={0 3mm 0 0},clip]{\currfiledir 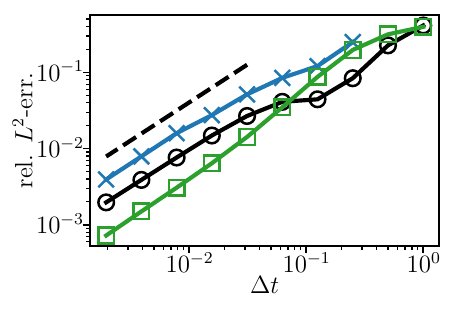}
		\caption{CRN model, $\Delta x=0.2 \mm$.}
		\label{fig:cub2d_conv_ref_0_CRN}
	\end{subfigure}\hfill%
	\begin{subfigure}[t]{\subfigsizet\textwidth}
		\centering
		\includegraphics[scale=\plotimscalet,trim={0 3mm 0 0},clip]{\currfiledir 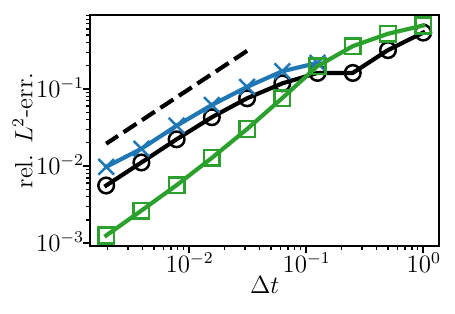}
		\caption{TTP model, $\Delta x=0.2 \mm$.}
		\label{fig:cub2d_conv_ref_0_TTP}
	\end{subfigure}
	\begin{subfigure}[t]{\subfigsizet\textwidth}
		\centering
		\includegraphics[scale=\plotimscalet,trim={0 3mm 0 0},clip]{\currfiledir 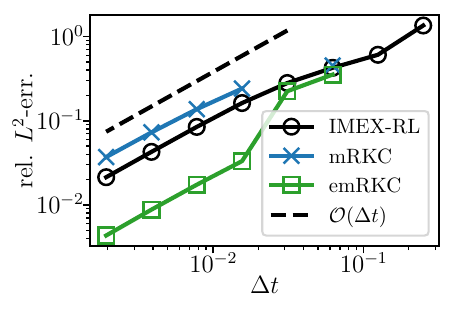}
		\caption{HH model, $\Delta x=0.1 \mm$.}
		\label{fig:cub2d_conv_ref_1_HH}
	\end{subfigure}\hfill%
	\begin{subfigure}[t]{\subfigsizet\textwidth}
		\centering
		\includegraphics[scale=\plotimscalet,trim={0 3mm 0 0},clip]{\currfiledir 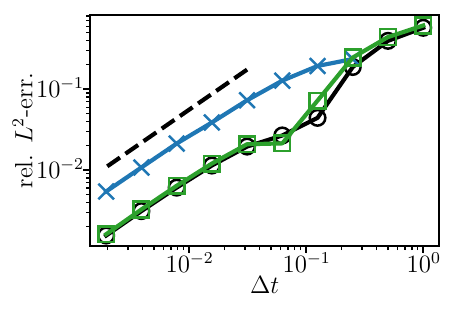}
		\caption{CRN model, $\Delta x=0.1 \mm$.}
		\label{fig:cub2d_conv_ref_1_CRN}
	\end{subfigure}\hfill%
	\begin{subfigure}[t]{\subfigsizet\textwidth}
		\centering
		\includegraphics[scale=\plotimscalet,trim={0 3mm 0 0},clip]{\currfiledir 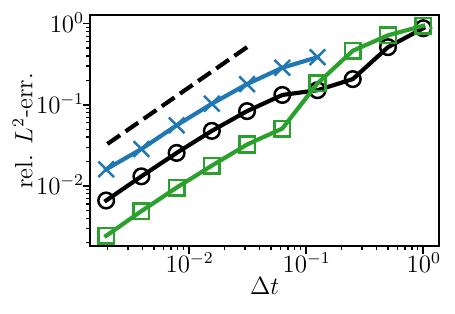}
		\caption{TTP model, $\Delta x=0.1 \mm$.}
		\label{fig:cub2d_conv_ref_1_TTP}
	\end{subfigure}
	\begin{subfigure}[t]{\subfigsizet\textwidth}
		\centering
		\includegraphics[scale=\plotimscalet,trim={0 3mm 0 0},clip]{\currfiledir 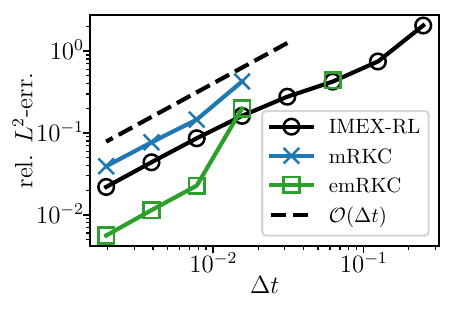}
		\caption{HH model, $\Delta x=0.05 \mm$.}
		\label{fig:cub2d_conv_ref_2_HH}
	\end{subfigure}\hfill%
	\begin{subfigure}[t]{\subfigsizet\textwidth}
		\centering
		\includegraphics[scale=\plotimscalet,trim={0 3mm 0 0},clip]{\currfiledir 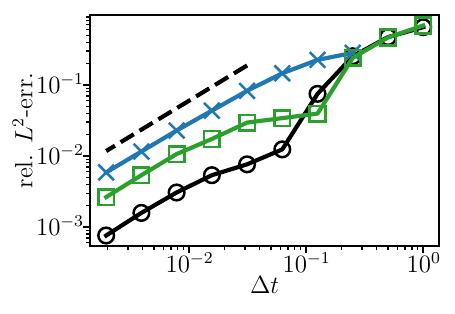}
		\caption{CRN model, $\Delta x=0.05 \mm$.}
		\label{fig:cub2d_conv_ref_2_CRN}
	\end{subfigure}\hfill%
		\begin{subfigure}[t]{\subfigsizet\textwidth}
		\centering
		\includegraphics[scale=\plotimscalet,trim={0 3mm 0 0},clip]{\currfiledir 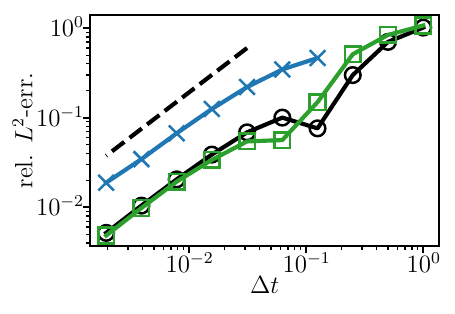}
		\caption{TTP model, $\Delta x=0.05 \mm$.}
		\label{fig:cub2d_conv_ref_2_TTP}
	\end{subfigure}
	\caption{Two dimensional rectangle. Convergence experiments with different mesh sizes $\Delta x$ and ionic models.}
	\label{fig:cub2d_conv}
\end{figure}

To begin we assess the accuracy of the various time integration methods, in particular we verify that the first-order of convergence predicted by \cref{thm:convergence} also holds in practice. More precisely, for each mesh size $\Delta x=\SI{0.2}{\mm}$, $\SI{0.1}{\mm}$, $\SI{0.05}{\mm}$ and every ionic model HH, CRN, TTP, we successively reduce the step size as $\Delta t_i=2^{-i}\si{\milli\second}$ for $i=0,\ldots,9$, i.e., from $\SI{1}{\milli\second}$ down to approximately $\SI{0.002}{\milli\second}$. As shown in \cref{fig:cub2d_conv}, all methods achieve the expected first-order convergence.
Moreover, in all but one case, the emRKC is equally or even more accurate than the standard IMEX-RL method.

%The convergence results for all mesh sizes and ionic models are collected in %\cref{fig:cub2d_conv}. We see that the expected first-order of convergence is %confirmed in all cases and that, except for one case, the emRKC is more precise or %equivalent to the baseline method IMEX-RL.

Nonetheless, in \cref{fig:cub2d_conv_ref_0_HH,fig:cub2d_conv_ref_1_HH,fig:cub2d_conv_ref_2_HH} we also notice that some data points for large $\Delta t$ are missing
for the mRKC and emRKC methods, indicating instability, yet only for the HH model and not for the more realistic CRN and TTP models. We believe that instabilities in the HH model are due to the overly large imaginary part of the eigenvalues of the right-hand side's Jacobian, compared to the more realistic CRN and TTP models whose eigenvalues have relatively small imaginary parts \cite{Spiteri2010,SpD12}.
In fact, those instabilities become even more pronounced as the mesh size decreases, 
probably because the stability of mRKC and emRKC for non-real eigenvalues deteriorates as $f_F$ becomes increasingly stiffer; this effect, however, is not well understood yet, since a stability analysis of mRKC and emRKC with complex eigenvalues is still lacking.

\begin{figure}
	\begin{subfigure}[t]{\subfigsizet\textwidth}
		\centering
		\includegraphics[scale=\plotimscalet,trim={0 3mm 0 0},clip]{\currfiledir 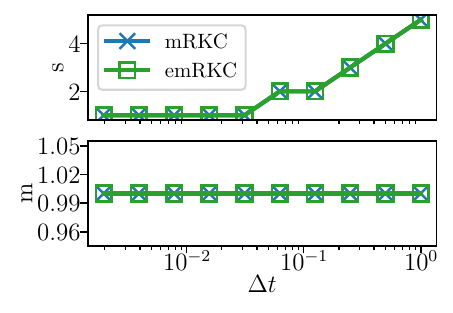}
		\caption{HH model, $\Delta x=0.2 \mm$.}
		\label{fig:cub2d_stages_ref_0_HH}
	\end{subfigure}\hfill%
	\begin{subfigure}[t]{\subfigsizet\textwidth}
		\centering
		\includegraphics[scale=\plotimscalet,trim={0 3mm 0 0},clip]{\currfiledir 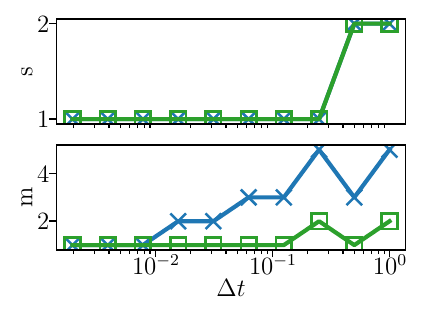}
		\caption{CRN model, $\Delta x=0.2 \mm$.}
		\label{fig:cub2d_stages_ref_0_CRN}
	\end{subfigure}\hfill%
	\begin{subfigure}[t]{\subfigsizet\textwidth}
		\centering
		\includegraphics[scale=\plotimscalet,trim={0 3mm 0 0},clip]{\currfiledir 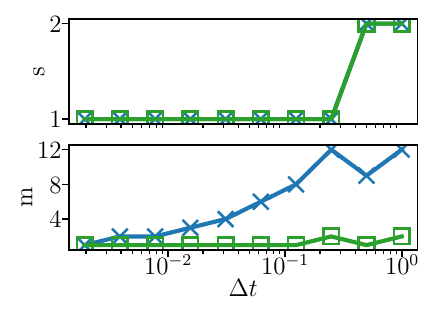}
		\caption{TTP model, $\Delta x=0.2 \mm$.}
		\label{fig:cub2d_stages_ref_0_TTP}
	\end{subfigure}\\
	\begin{subfigure}[t]{\subfigsizet\textwidth}
		\centering
		\includegraphics[scale=\plotimscalet,trim={0 3mm 0 0},clip]{\currfiledir 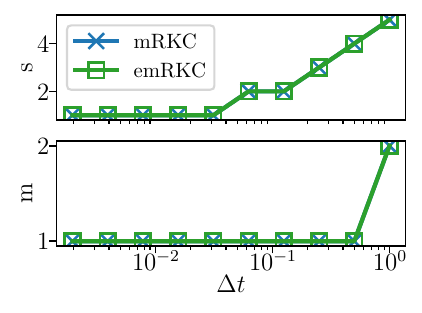}
		\caption{HH model, $\Delta x=0.1 \mm$.}
		\label{fig:cub2d_stages_ref_1_HH}
	\end{subfigure}\hfill%
	\begin{subfigure}[t]{\subfigsizet\textwidth}
		\centering
		\includegraphics[scale=\plotimscalet,trim={0 3mm 0 0},clip]{\currfiledir 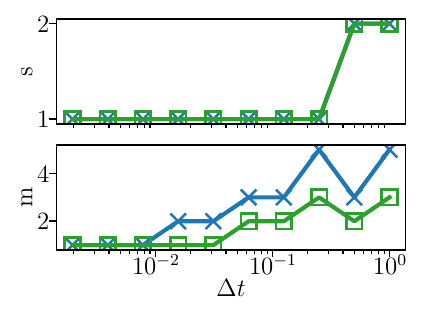}
		\caption{CRN model, $\Delta x=0.1 \mm$.}
		\label{fig:cub2d_stages_ref_1_CRN}
	\end{subfigure}\hfill%
	\begin{subfigure}[t]{\subfigsizet\textwidth}
		\centering
		\includegraphics[scale=\plotimscalet,trim={0 3mm 0 0},clip]{\currfiledir 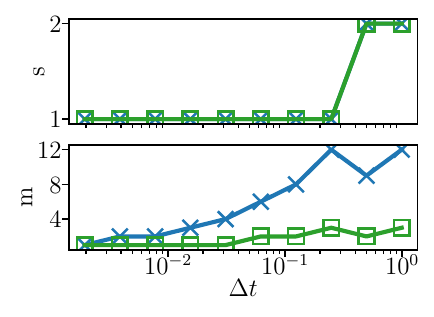}
		\caption{TTP model, $\Delta x=0.1 \mm$.}
		\label{fig:cub2d_stages_ref_1_TTP}
	\end{subfigure}\\
	\begin{subfigure}[t]{\subfigsizet\textwidth}
		\centering
		\includegraphics[scale=\plotimscalet,trim={0 3mm 0 0},clip]{\currfiledir 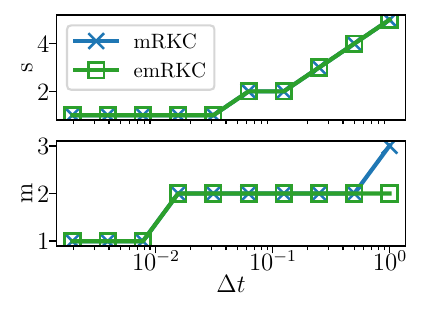}
		\caption{HH model, $\Delta x=0.05 \mm$.}
		\label{fig:cub2d_stages_ref_2_HH}
	\end{subfigure}\hfill%
	\begin{subfigure}[t]{\subfigsizet\textwidth}
		\centering
		\includegraphics[scale=\plotimscalet,trim={0 3mm 0 0},clip]{\currfiledir 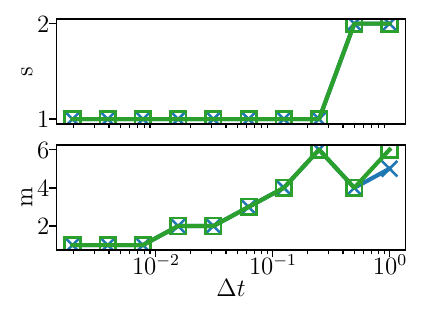}
		\caption{CRN model, $\Delta x=0.05 \mm$.}
		\label{fig:cub2d_stages_ref_2_CRN}
	\end{subfigure}\hfill%
	\begin{subfigure}[t]{\subfigsizet\textwidth}
		\centering
		\includegraphics[scale=\plotimscalet,trim={0 3mm 0 0},clip]{\currfiledir 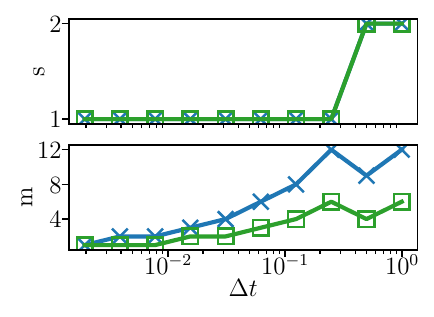}
		\caption{TTP model, $\Delta x=0.05 \mm$.}
		\label{fig:cub2d_stages_ref_2_TTP}
	\end{subfigure}
	\caption{Two dimensional rectangle. Number of stages taken by the mRKC and emRKC methods, with different mesh sizes $\Delta x$ and ionic models. Recall that $s$, $m$ are the number of stages taken by the outer and inner RKC iterations, respectively.}
	\label{fig:cub2d_stages}
\end{figure}

In \cref{fig:cub2d_stages}, we display the number of stages $s,m$ used by 
the mRKC and emRKC methods.
For the CRN and TTP models, a closer look reveals that both tend to lose some accuracy when the number of inner stages $m$ increases. As expected, the number of outer stages $s$ typically remains the same for both methods, since
$s$ depends on $f_S$ for emRKC but on $f_S+f_E^2$ for mRKC, yet with $f_E^2$ much less stiff than $f_S$. The number of inner stages $m$ for emRKC depends only on the diffusion term $f_F$, whereas for mRKC it depends on $f_F+f_E^1$. In \cref{fig:cub2d_stages}, we also observe that the number of stages $m$ for the two methods often differ, thus indicating that $f_E^1$ is indeed stiffer than the diffusion $f_F$. In those instances, the efficiency of mRKC deteriorates against emRKC. We also observe for mRKC that $m$ depends not only on the ionic model but also on the mesh size, while for emRKC it depends only on the mesh size. Finally, the ``V'' shape behavior 
in the graph of $m$ apparent at larger step sizes stems from the ``multiplicative'' stabilization procedure of the methods, as discussed in \cref{sec:mES}. Indeed, for decreasing $\Delta t$, a constant number of outer stages $s$ also implies 
a decreasing number of inner stages $m$. In contrast, for decreasing $s$, the number of inner stages $m$ can slightly increase to compensate for the weaker stabilization of the outer loop. Still, $s \times m$ always decreases with $\Delta t$. 

\begin{figure}
	\begin{subfigure}[t]{\subfigsizet\textwidth}
		\centering
		\includegraphics[scale=\plotimscalet,trim={0 3mm 0 0},clip]{\currfiledir 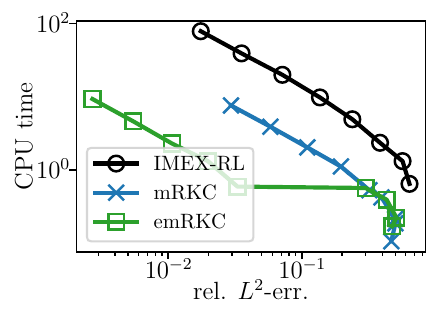}
		\caption{HH model, $\Delta x=0.2 \mm$.}
		\label{fig:cub2d_eff_ref_0_HH}
	\end{subfigure}\hfill%
	\begin{subfigure}[t]{\subfigsizet\textwidth}
		\centering
		\includegraphics[scale=\plotimscalet,trim={0 3mm 0 0},clip]{\currfiledir 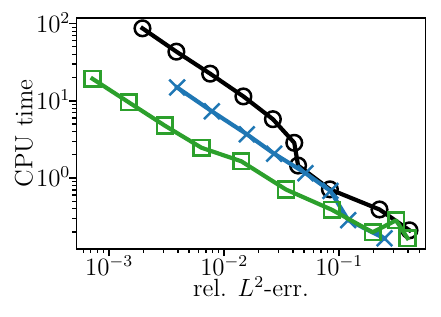}
		\caption{CRN model, $\Delta x=0.2 \mm$.}
		\label{fig:cub2d_eff_ref_0_CRN}
	\end{subfigure}\hfill%
	\begin{subfigure}[t]{\subfigsizet\textwidth}
		\centering
		\includegraphics[scale=\plotimscalet,trim={0 3mm 0 0},clip]{\currfiledir 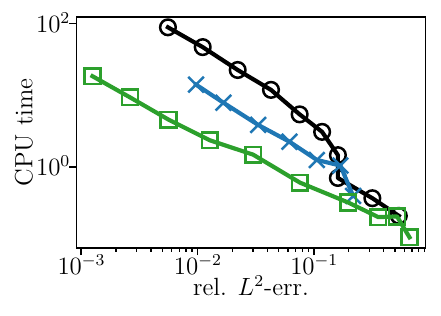}
		\caption{TTP model, $\Delta x=0.2 \mm$.}
		\label{fig:cub2d_eff_ref_0_TTP}
	\end{subfigure}
	\begin{subfigure}[t]{\subfigsizet\textwidth}
		\centering
		\includegraphics[scale=\plotimscalet,trim={0 3mm 0 0},clip]{\currfiledir 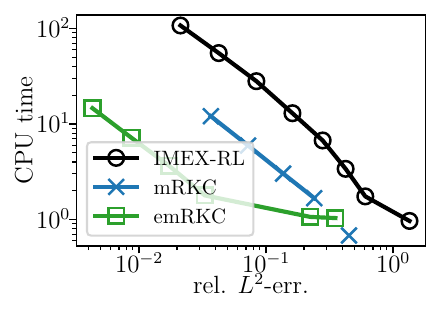}
		\caption{HH model, $\Delta x=0.1 \mm$.}
		\label{fig:cub2d_eff_ref_1_HH}
	\end{subfigure}\hfill%
	\begin{subfigure}[t]{\subfigsizet\textwidth}
		\centering
		\includegraphics[scale=\plotimscalet,trim={0 3mm 0 0},clip]{\currfiledir 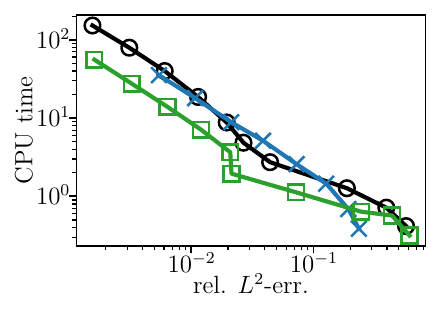}
		\caption{CRN model, $\Delta x=0.1 \mm$.}
		\label{fig:cub2d_eff_ref_1_CRN}
	\end{subfigure}\hfill%
	\begin{subfigure}[t]{\subfigsizet\textwidth}
		\centering
		\includegraphics[scale=\plotimscalet,trim={0 3mm 0 0},clip]{\currfiledir 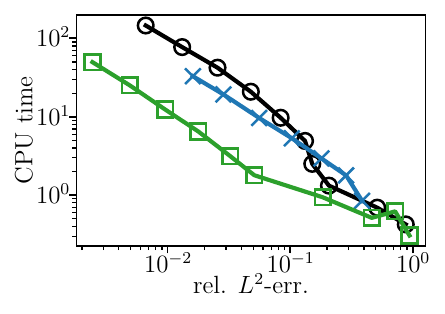}
		\caption{TTP model, $\Delta x=0.1 \mm$.}
		\label{fig:cub2d_eff_ref_1_TTP}
	\end{subfigure}
	\begin{subfigure}[t]{\subfigsizet\textwidth}
		\centering
		\includegraphics[scale=\plotimscalet,trim={0 3mm 0 0},clip]{\currfiledir 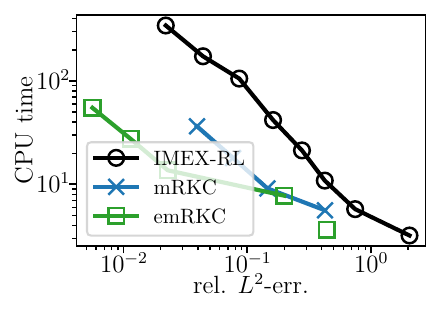}
		\caption{HH model, $\Delta x=0.05 \mm$.}
		\label{fig:cub2d_eff_ref_2_HH}
	\end{subfigure}\hfill%
	\begin{subfigure}[t]{\subfigsizet\textwidth}
		\centering
		\includegraphics[scale=\plotimscalet,trim={0 3mm 0 0},clip]{\currfiledir 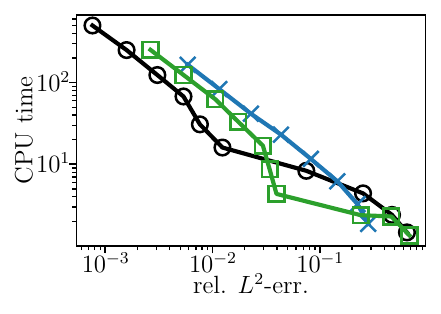}
		\caption{CRN model, $\Delta x=0.05 \mm$.}
		\label{fig:cub2d_eff_ref_2_CRN}
	\end{subfigure}\hfill%
	\begin{subfigure}[t]{\subfigsizet\textwidth}
		\centering
		\includegraphics[scale=\plotimscalet,trim={0 3mm 0 0},clip]{\currfiledir 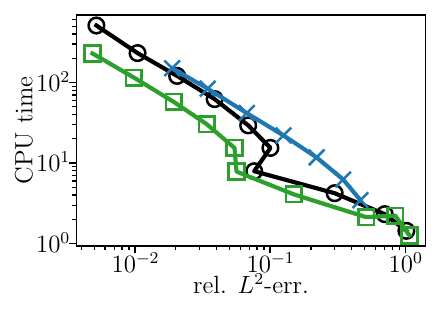}
		\caption{TTP model, $\Delta x=0.05 \mm$.}
		\label{fig:cub2d_eff_ref_2_TTP}
	\end{subfigure}
	\caption{Two dimensional rectangle. Efficiency experiments with different mesh sizes $\Delta x$ and ionic models.}
	\label{fig:cub2d_eff}
\end{figure}

Finally, in \cref{fig:cub2d_eff} we display work vs. accuracy by monitoring the $L^2$-norm relative errors against the total CPU time. The emRKC method tends to be the most efficient scheme except for one instance, while the mRKC method also performs reasonably well. The high performance of emRKC compared to the other methods stems from a number of factors. First, in comparison to IMEX-RL, the explicit emRKC method is
cheaper and tends to be more accurate because the $f_S, f_E$ terms take part in the stabilization procedure for the diffusion term $f_F$ (\cref{algo:emRKC_inner} of \cref{algo:emRKC}); thus, a strong coupling between the different model's components is guaranteed (not so in the splitting strategy of IMEX-RL). Second, in comparison to mRKC, emRKC is cheaper due to the smaller number of inner stages $m$ and has higher accuracy because it employs the more accurate exponential integrator for the gating variables. 

\subsection{Three-dimensional cuboid}\label{sec:cuboid3d}
Next, we perform a series of numerical experiments very similar to those in \cref{sec:cuboid2d} but in three dimensions. Hence, we set $\Omega=[0,20]\times [0,7]\times [0,3]\,\si{\milli\metre\cubed}$ while the applied stimulus and final simulation time remain identical.
Again all methods achieved first-order convergence in time for all three ionic models -- those plots are omitted here for the sake of brevity. We also omit the graphs of the number of stages taken by mRKC and emRKC, which are identical to the two-dimensional case. 

\begin{figure}
	\begin{subfigure}[t]{\subfigsizet\textwidth}
		\centering
		\includegraphics[scale=\plotimscalet,trim={0 3mm 0 0},clip]{\currfiledir 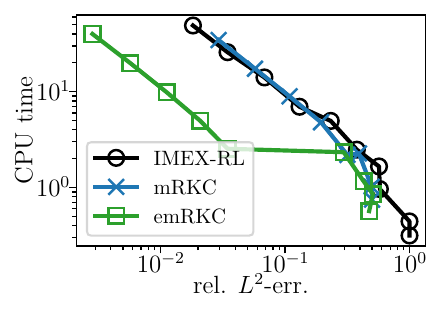}
		\caption{HH model, $\Delta x=0.2 \mm$.}
		\label{fig:cub3d_eff_ref_0_HH}
	\end{subfigure}\hfill%
	\begin{subfigure}[t]{\subfigsizet\textwidth}
		\centering
		\includegraphics[scale=\plotimscalet,trim={0 3mm 0 0},clip]{\currfiledir 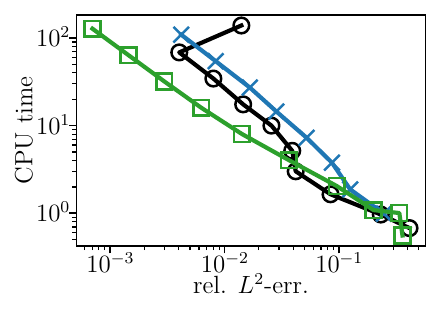}
		\caption{CRN model, $\Delta x=0.2 \mm$.}
		\label{fig:cub3d_eff_ref_0_CRN}
	\end{subfigure}\hfill%
	\begin{subfigure}[t]{\subfigsizet\textwidth}
		\centering
		\includegraphics[scale=\plotimscalet,trim={0 3mm 0 0},clip]{\currfiledir 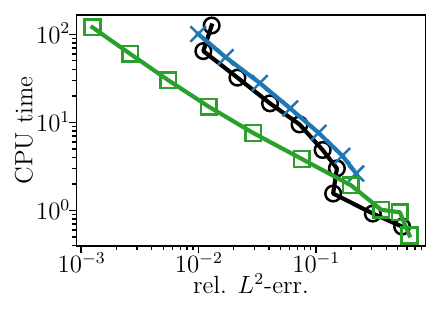}
		\caption{TTP model, $\Delta x=0.2 \mm$.}
		\label{fig:cub3d_eff_ref_0_TTP}
	\end{subfigure}
	\begin{subfigure}[t]{\subfigsizet\textwidth}
		\centering
		\includegraphics[scale=\plotimscalet,trim={0 3mm 0 0},clip]{\currfiledir 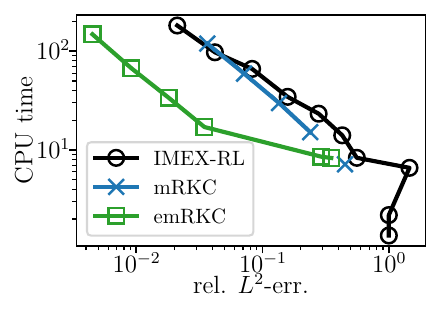}
		\caption{HH model, $\Delta x=0.1 \mm$.}
		\label{fig:cub3d_eff_ref_1_HH}
	\end{subfigure}\hfill%
	\begin{subfigure}[t]{\subfigsizet\textwidth}
		\centering
		\includegraphics[scale=\plotimscalet,trim={0 3mm 0 0},clip]{\currfiledir 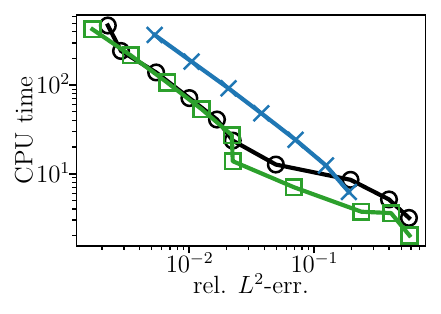}
		\caption{CRN model, $\Delta x=0.1 \mm$.}
		\label{fig:cub3d_eff_ref_1_CRN}
	\end{subfigure}\hfill%
	\begin{subfigure}[t]{\subfigsizet\textwidth}
		\centering
		\includegraphics[scale=\plotimscalet,trim={0 3mm 0 0},clip]{\currfiledir 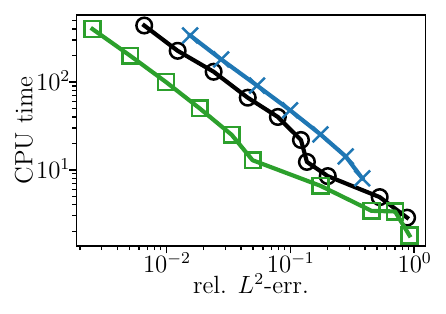}
		\caption{TTP model, $\Delta x=0.1 \mm$.}
		\label{fig:cub3d_eff_ref_1_TTP}
	\end{subfigure}
	\begin{subfigure}[t]{\subfigsizet\textwidth}
		\centering
		\includegraphics[scale=\plotimscalet,trim={0 3mm 0 0},clip]{\currfiledir 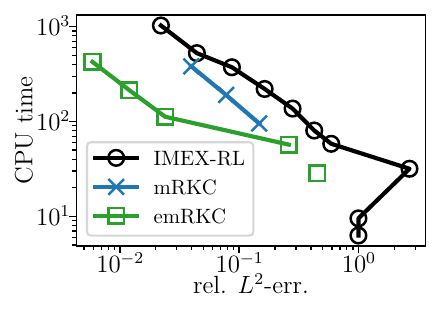}
		\caption{HH model, $\Delta x=0.05 \mm$.}
		\label{fig:cub3d_eff_ref_2_HH}
	\end{subfigure}\hfill%
	\begin{subfigure}[t]{\subfigsizet\textwidth}
		\centering
		\includegraphics[scale=\plotimscalet,trim={0 3mm 0 0},clip]{\currfiledir 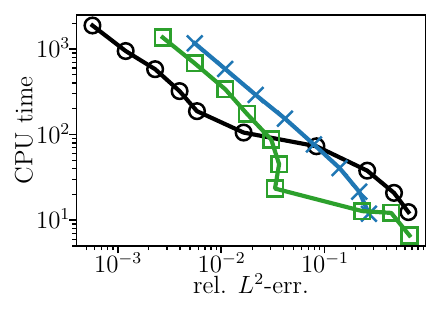}
		\caption{CRN model, $\Delta x=0.05 \mm$.}
		\label{fig:cub3d_eff_ref_2_CRN}
	\end{subfigure}\hfill%
	\begin{subfigure}[t]{\subfigsizet\textwidth}
		\centering
		\includegraphics[scale=\plotimscalet,trim={0 3mm 0 0},clip]{\currfiledir 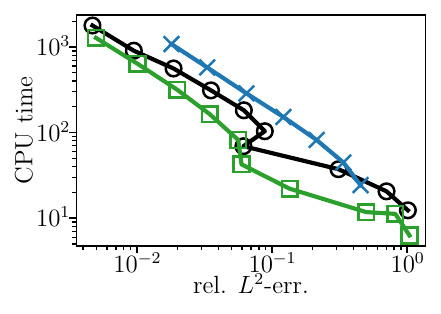}
		\caption{TTP model, $\Delta x=0.05 \mm$.}
		\label{fig:cub3d_eff_ref_2_TTP}
	\end{subfigure}
	\caption{Three dimensional cuboid. Efficiency experiments with different mesh sizes $\Delta x$ and ionic models.}
	\label{fig:cub3d_eff}
\end{figure}

The work vs. accuracy diagrams, shown in \cref{fig:cub3d_eff}, again confirm the high efficiency of emRKC. Still, it tends to be slower than in 2-D probably because the inner iterations are more expensive when the diffusion matrix stems from a three-dimensional problem. Again we note that the mRKC and emRKC schemes are stable for the CRN and TTP models, even at the largest step size $\Delta t=\SI{1}{\milli\second}$, but not for the HH model.

\subsection{Left atrium geometry}\label{sec:fastl_LA}
\begin{figure}
    \begin{subfigure}[t]{0.09\textwidth}
    \centering
    \includegraphics[width=\textwidth]{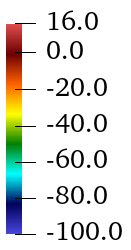}    
    \end{subfigure}\hfill%
    \begin{subfigure}[t]{0.44\textwidth}
    \centering
    \includegraphics[width=0.5\textwidth]{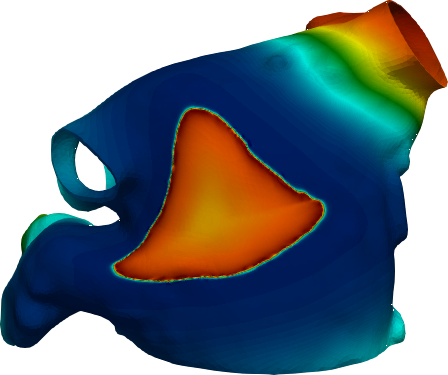}%
    \includegraphics[width=0.5\textwidth]{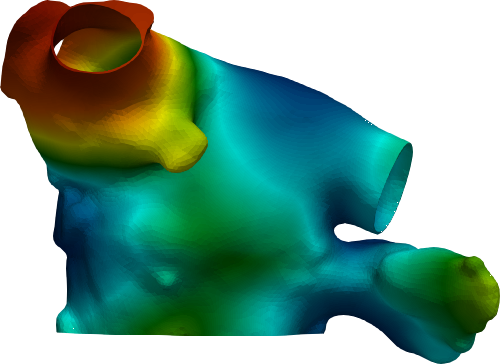}
    \caption{Initial value.}
    \label{fig:init_val}
    \end{subfigure}\hfill%
    \begin{subfigure}[t]{0.44\textwidth}
    \centering
    \includegraphics[width=0.5\textwidth]{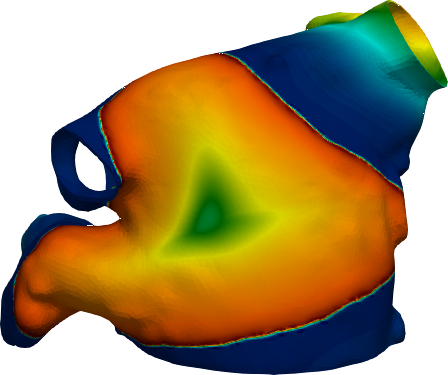}%
    \includegraphics[width=0.5\textwidth]{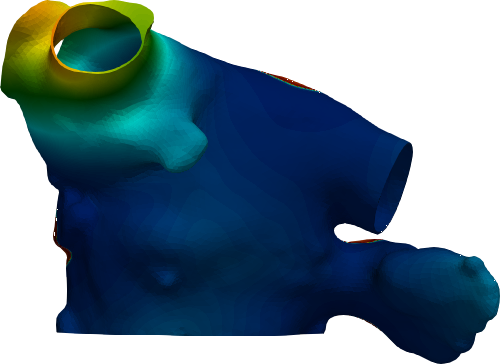}
    \caption{Solution at $t=\SI{50}{\milli\second}$.}
    \label{fig:end_val}
    \end{subfigure}\hfill%
    \caption{Left atrium geometry. Initial value and solution at $t=\SI{50}{\milli\second}$ for the CRN ionic model.}
\end{figure}

Here we compare the performance of the mRKC and emRKC methods against the IMEX-RL benchmark but now consider a realistic left atrium geometry. The computational mesh was obtained by
refining twice the realistic mesh used in \cite{FasTobCroWhiRajMcCSanHoONPlaBisNie18,Fastl2021}, which led to an average mesh size of $\SI{0.26}{\mm}$, maximal mesh size of $\SI{0.69}{\mm}$, minimal mesh of size $\SI{0.026}{\mm}$, and about 1.7 million degrees of freedom for the FE-solution.
The vector fields $\bm a_l(\bm x)$ and $\bm a_t(\bm x)$ which define the fiber and sheet orientation in \cref{eq:def_D}, have been generated with the lifex-fiber library \cite{AfrPieFedDedQua23,piersanti2021modeling}. The initial value, 
shown in \cref{fig:init_val}, is no longer constant throughout $\Omega$ but instead generated by solving the monodomain equation for $t\in[0,800]\,\si{\milli\second}$ and applying the same stimulus-pacing protocol as in \cite{GanPezGhaKraPerSah22,CostabalFKS23} in a fixed region of the domain. 

\begin{figure}	
	\begin{subfigure}[t]{\subfigsizet\textwidth}
		\centering
		\includegraphics[scale=\plotimscalet,trim={0 3mm 0 0},clip]{\currfiledir 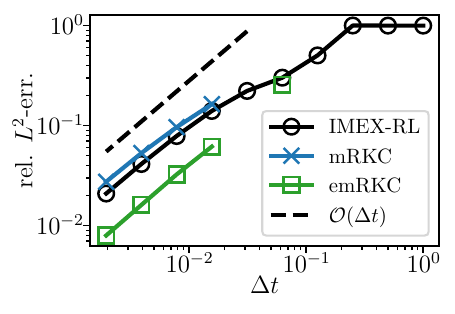}
		\caption{HH model.}
		\label{fig:03_fastl_LA_conv_ref_2_HH}
	\end{subfigure}\hfill%
	\begin{subfigure}[t]{\subfigsizet\textwidth}
		\centering
		\includegraphics[scale=\plotimscalet,trim={0 3mm 0 0},clip]{\currfiledir 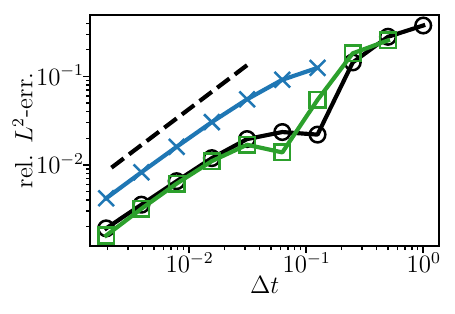}
		\caption{CRN model.}
		\label{fig:03_fastl_LA_conv_ref_2_CRN}
	\end{subfigure}\hfill%
		\begin{subfigure}[t]{\subfigsizet\textwidth}
		\centering
		\includegraphics[scale=\plotimscalet,trim={0 3mm 0 0},clip]{\currfiledir 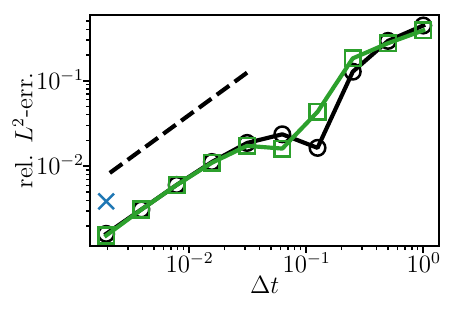}
		\caption{TTP model.}
		\label{fig:03_fastl_LA_conv_ref_2_TTP}
	\end{subfigure}
	\caption{Left atrium geometry. Convergence experiments with different ionic models.}
	\label{fig:03_fastl_LA_conv}
\end{figure}

As in \cref{sec:cuboid2d,sec:cuboid3d}, we perform convergence and efficiency experiments for the three ionic models HH, CRN, and TTP. The final time is set to $T=\SI{50}{\milli\second}$, see \cref{fig:end_val} for an illustration of the solution at $t=T$. Here, the applied stimulus is identically zero, as the initial condition already contains a propagating wave. Again, all three methods achieve first-order convergence in time for all three models, as shown
in \cref{fig:03_fastl_LA_conv}. Unlike the previous cases, however, the mRKC method here is unstable for the TTP model, whereas the emRKC method always remains stable.
The relative error for emRKC is always comparable if not better (HH case) to the error for IMEX-RL, which is the standard method in the cardiac community. Note that relative error of 1\% to 5\% is comparable to the spatial error, which is of the same order for linear FE and mesh size of \SI{\approx 0.25}{\milli\meter}, see~\cite{PHS16}.

\begin{figure}	
	\begin{subfigure}[t]{\subfigsizet\textwidth}
		\centering
		\includegraphics[scale=\plotimscalet,trim={0 3mm 0 0},clip]{\currfiledir 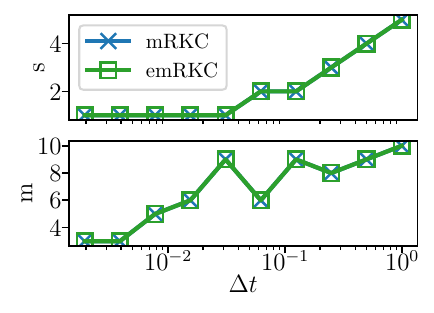}
		\caption{HH model.}
		\label{fig:03_fastl_LA_stages_ref_2_HH}
	\end{subfigure}\hfill%
	\begin{subfigure}[t]{\subfigsizet\textwidth}
		\centering
		\includegraphics[scale=\plotimscalet,trim={0 3mm 0 0},clip]{\currfiledir 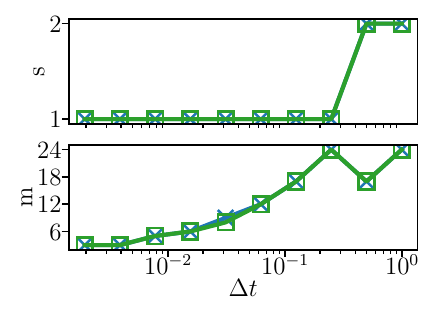}
		\caption{CRN model.}
		\label{fig:03_fastl_LA_stages_ref_2_CRN}
	\end{subfigure}\hfill%
	\begin{subfigure}[t]{\subfigsizet\textwidth}
		\centering
		\includegraphics[scale=\plotimscalet,trim={0 3mm 0 0},clip]{\currfiledir 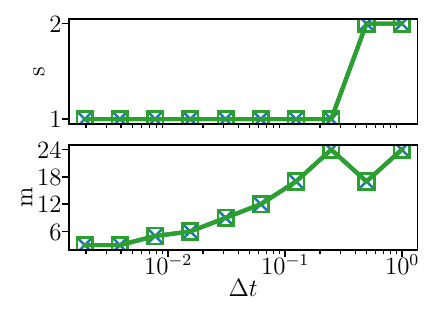}
		\caption{TTP model.}
		\label{fig:03_fastl_LA_stages_ref_2_TTP}
	\end{subfigure}
	\caption{Left atrium geometry. Number of stages used by the mRKC and emRKC methods with different ionic models. Recall that $s$, $m$ are the number of stages taken by the outer and inner RKC iterations, respectively.}
	\label{fig:03_fastl_LA_stages}
\end{figure}

\begin{figure}	
	\begin{subfigure}[t]{\subfigsizet\textwidth}
		\centering
		\includegraphics[scale=\plotimscalet,trim={0 3mm 0 0},clip]{\currfiledir 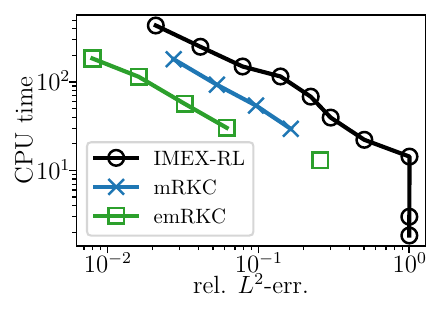}
		\caption{HH model.}
		\label{fig:03_fastl_LA_eff_ref_2_HH}
	\end{subfigure}\hfill%
	\begin{subfigure}[t]{\subfigsizet\textwidth}
		\centering
		\includegraphics[scale=\plotimscalet,trim={0 3mm 0 0},clip]{\currfiledir 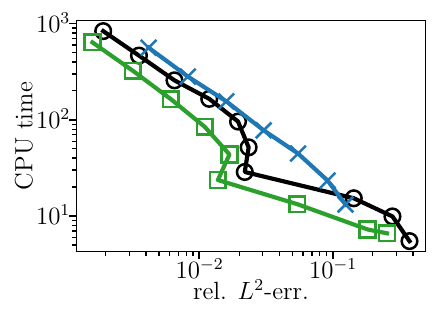}
		\caption{CRN model.}
		\label{fig:03_fastl_LA_eff_ref_2_CRN}
	\end{subfigure}\hfill%
	\begin{subfigure}[t]{\subfigsizet\textwidth}
		\centering
		\includegraphics[scale=\plotimscalet,trim={0 3mm 0 0},clip]{\currfiledir 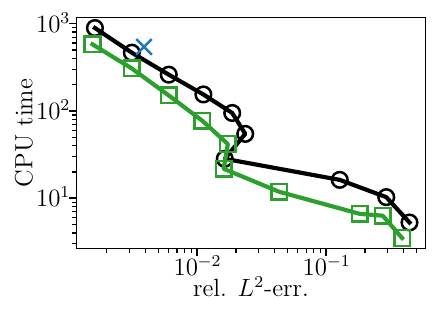}
		\caption{TTP model.}
		\label{fig:03_fastl_LA_eff_ref_2_TTP}
	\end{subfigure}
	\caption{Left atrium geometry. Efficiency experiments with different ionic models.}
	\label{fig:03_fastl_LA_eff}
\end{figure}

In \cref{fig:03_fastl_LA_stages}, we observe that the number of stages used by the mRKC and emRKC methods is higher because of the smaller minimal mesh size. Therefore, the mRKC method evaluates the stiff part $f_E^1$ of the ionic model many times \cref{eq:mRCK_fFfS_monodomain}, which may be the source of internal instabilities. As the internal stability analysis of mRKC is not available, this issue is not fully understood yet. Finally, we note that both explicit multirate methods always take a number of stages significantly larger than one, indicating 
that any standard explicit method would be highly unstable for such large step sizes. 
In \cref{fig:03_fastl_LA_eff}, the results from the efficiency experiment again suggest that emRKC is slightly faster than the IMEX-RL method, despite the high number of stages (\cref{fig:03_fastl_LA_stages}) due to the mesh induced
increased stiffness. 

\begin{figure}	
        \centering
        \begin{subfigure}[t]{\subfigsizet\textwidth}
		\centering
		\includegraphics[scale=\plotimscalet,trim={0mm 5mm 2mm 3mm},clip]{\currfiledir 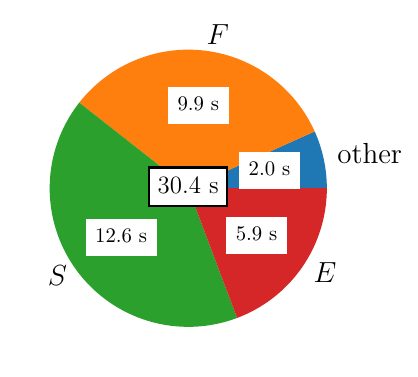}
		\caption{HH model.}
		\label{fig:prof_03_fastl_HH_emRKC}
	\end{subfigure}\hfill%	
        \begin{subfigure}[t]{\subfigsizet\textwidth}
		\centering
		\includegraphics[scale=\plotimscalet,trim={0mm 5mm 2mm 3mm},clip]{\currfiledir 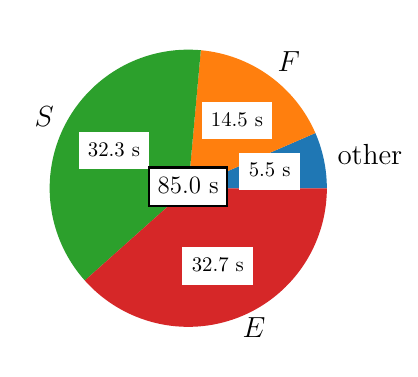}
		\caption{CRN model.}
		\label{fig:prof_03_fastl_CRN_emRKC}
	\end{subfigure}\hfill%
        \begin{subfigure}[t]{\subfigsizet\textwidth}
		\centering
		\includegraphics[scale=\plotimscalet,trim={0mm 5mm 2mm 3mm},clip]{\currfiledir 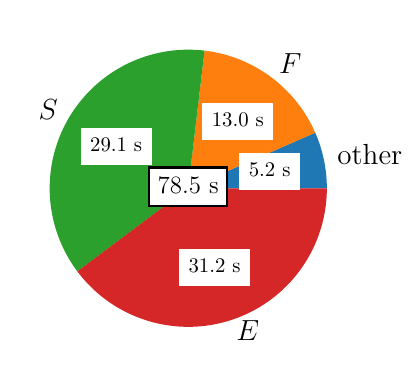}
		\caption{TTP model.}
		\label{fig:prof_03_fastl_TTP_emRKC}
	\end{subfigure}\hfill
 	\caption{Left atrium geometry. Code profiling of the emRKC method. Wedges $F$, $S$, $E$ indicate time spent evaluating $f_F$, $f_S$, and computing the exponential step, respectively. Remaining operations are represented by the ``other'' wedge.}
	\label{fig:prof_03_fastl_emRKC}
\end{figure}

Finally, we provide some code profiling results for the emRKC and IMEX-RL methods. First, we measure how much time the emRKC method spends in evaluating $f_F$, $f_S$, and in performing the exponential Euler step for $f_E$ (see \cref{algo:emRKC}). To do so, we apply emRKC with $\Delta t=2^{-6}\,\si{\milli\second}=\SI{0.015625}{\milli\second}$ and profile the \textsc{emRKC\_Step} function of \cref{algo:emRKC}. We choose this $\Delta t$ since it is the largest one yielding a relative error of about 1\% (\cref{fig:03_fastl_LA_conv}). Moreover, in view of the efficiency results depicted in \cref{fig:03_fastl_LA_eff}, we deduce that the outcome of the profiling is only weakly affected by the size of $\Delta t$. The results for the three ionic models are displayed in \cref{fig:prof_03_fastl_emRKC}, where for each graph we show the average over five runs. The four wedges in the pie chart are as follows: $F$ represents the time spent in evaluating $f_F$, $S$ the time spent evaluating $f_S$, $E$ the time for computing the exponential step for $f_E$, and ``other'' the time required for all remaining operations, such as communications and vector additions. The size of each wedge is proportional to the time spent in the corresponding task, with the total time (spent in the \textsc{emRKC\_Step} routine) at the center of the pie chart. At first sight, the time spent in $F$ seems
to depend on the ionic model, even though the number of matrix-vector multiplications only depends on the mesh. Still, a closer look at those timings reveals that this variability is due to large fluctuations in the time spent computing matrix-vector multiplications.
In contrast, the execution timings for $S$ and $E$ are truly ionic model dependent. 

\begin{figure}	
        \centering
        \begin{subfigure}[t]{\subfigsizet\textwidth}
		\centering
		\includegraphics[scale=\plotimscalet,trim={0mm 5mm 2mm 3mm},clip]{\currfiledir 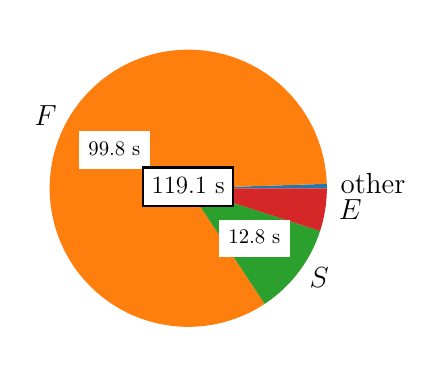}
		\caption{HH model.}
		\label{fig:prof_03_fastl_HH_IMEXRL}
	\end{subfigure}\hfill            
        \begin{subfigure}[t]{\subfigsizet\textwidth}
		\centering
		\includegraphics[scale=\plotimscalet,trim={0mm 5mm 2mm 3mm},clip]{\currfiledir 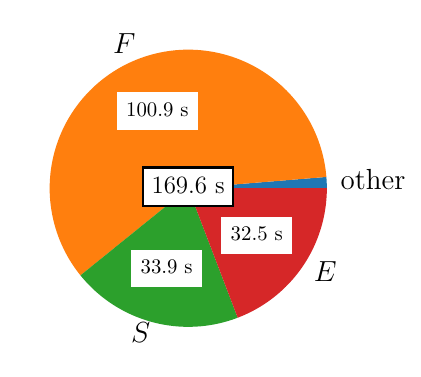}
		\caption{CRN model.}
		\label{fig:prof_03_fastl_CRN_IMEXRL}
	\end{subfigure}\hfill         
        \begin{subfigure}[t]{\subfigsizet\textwidth}
		\centering
		\includegraphics[scale=\plotimscalet,trim={0mm 5mm 2mm 3mm},clip]{\currfiledir 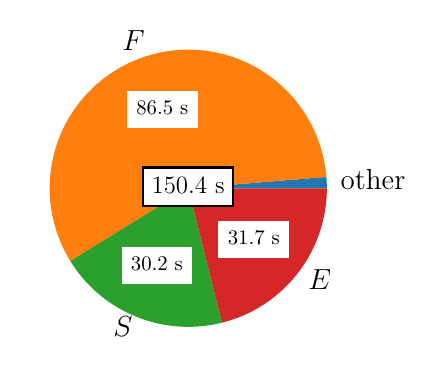}
		\caption{TTP model.}
		\label{fig:prof_03_fastl_TTP_IMEXRL}
	\end{subfigure}\hfill%	
	\caption{Left atrium geometry. Code profiling of the IMEX-RL method. Wedges $F$, $S$, $E$ indicate time spent solving the linear system for $f_F$, evaluating $f_S$, and computing the exponential step, respectively. Remaining operations are represented by the ``other'' wedge.}
	\label{fig:prof_03_fastl_IMEXRL}
\end{figure}

Next, we perform the same profiling experiments for the IMEX-RL code while measuring the time spent in \textsc{IMEX-RL\_Step} of \cref{algo:imexexp}. The corresponding pie charts are displayed in \cref{fig:prof_03_fastl_IMEXRL}, where the total execution time is again distributed between the four wedges $F$, $S$, $E$ and ``other''. Here, $F$ represents the time spent in solving the linear system in \cref{algo:imexexp} \cref{algo:imexexp_implicit}, $S$ the time spent in evaluating the slow components of the ionic model in \cref{algo:imexexp_slow} of \cref{algo:imexexp} , and $E$ the time spent in evaluating the exponential for the gating variables in \cref{algo:imexexp_exp} of \cref{algo:imexexp} . 

As expected, both emRKC and IMEX-RL spend approximately the same amount of time in $S$ and $E$, since the required operations are indeed identical. Due to the larger number of vector additions, the emRKC method spends more time in the ``other'' wedge. However, the total time required by IMEX-RL to solve the linear systems is about 6-7 times larger than the time spent by emRKC in evaluating $f_F$ (compare the two $F$ wedges). This significant reduction in execution time implies that the stabilization procedure of emRKC for the diffusion term $f_F$ is much cheaper than solving a linear system during each implicit Euler step. Remarkably, we achieve this reduction in computational cost without sacrificing accuracy, as illustrated in \cref{fig:03_fastl_LA_conv}. Overall, emRKC achieves a speed-up factor of 
about 2 over IMEX-RL
for the two most relevant ionic models CRN and TTP.

\subsection{Left atrium with fibrosis}\label{sec:fastl_LA_fibrosis}
\begin{figure}
    \centering
    \begin{subfigure}[t]{0.08\textwidth}
    \centering
    \includegraphics[width=\textwidth]{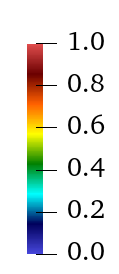}    
    \end{subfigure}
    \begin{subfigure}[t]{0.44\textwidth}
    \centering
    \includegraphics[width=0.5\textwidth]{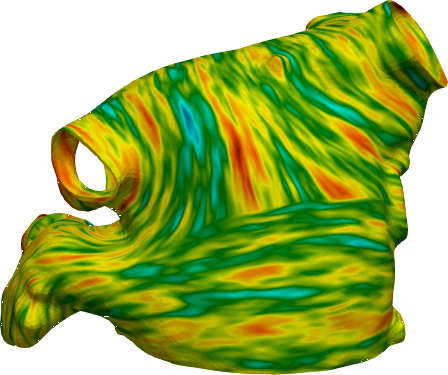}%
    \includegraphics[width=0.5\textwidth]{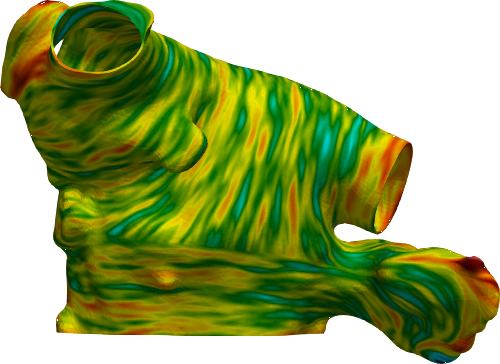}   
    \caption{Fibrosis random field $u(\bm x)$.}
    \label{fig:fibrosis_field}
    \end{subfigure}    
    \begin{subfigure}[t]{0.44\textwidth}
    \centering
    \includegraphics[width=0.5\textwidth]{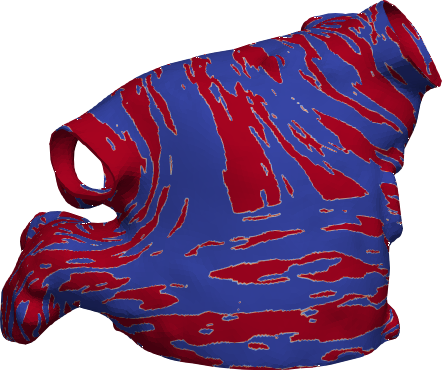}%
    \includegraphics[width=0.5\textwidth]{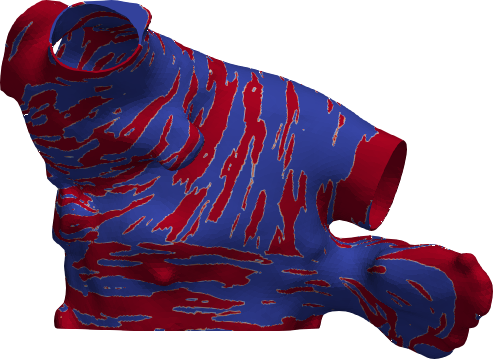} 
    \caption{Fibrosis pattern, in red where $\sigma_t=0$.}
    \label{fig:fibrosis_pattern}
    \end{subfigure}    
%    \caption{Left atrium with fibrosis. Fibrosis random field and pattern.}
%    \label{fig:fibrosis_field_pattern}
	\begin{subfigure}[t]{0.09\textwidth}
		\centering
		\includegraphics[width=\textwidth]{images/solutions/03_fastl_LA_bar.png}    
	\end{subfigure}\hfill%
	\begin{subfigure}[t]{0.44\textwidth}
		\centering
		\includegraphics[width=0.5\textwidth]{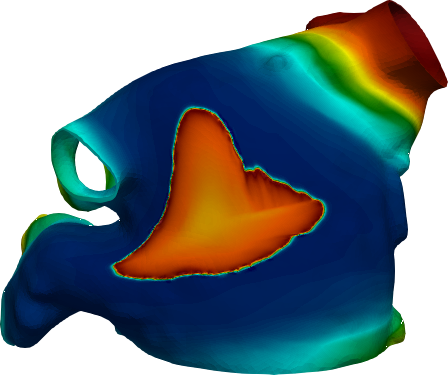}%
		\includegraphics[width=0.5\textwidth]{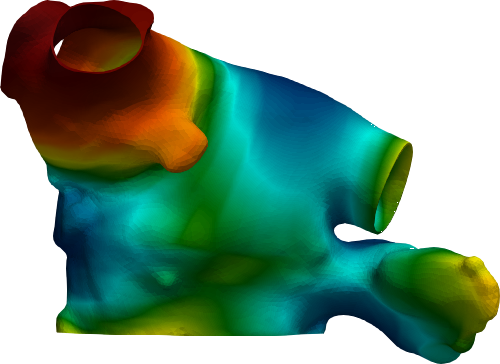}
		\caption{Initial value.}
		\label{fig:init_val_fibrosis}
	\end{subfigure}\hfill%
	\begin{subfigure}[t]{0.44\textwidth}
		\centering
		\includegraphics[width=0.5\textwidth]{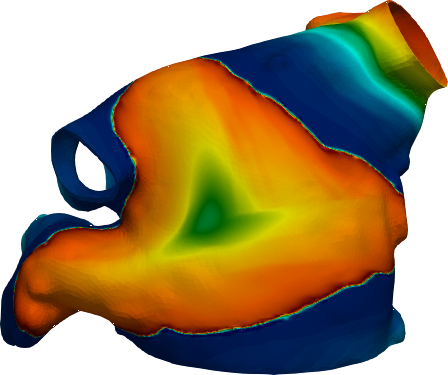}%
		\includegraphics[width=0.5\textwidth]{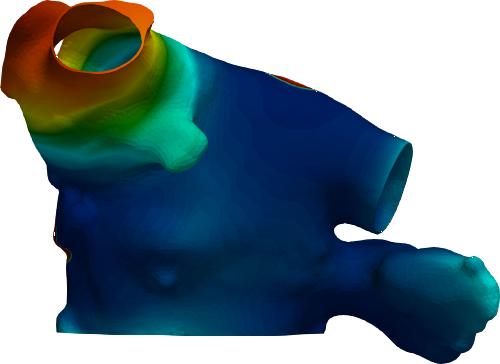}
		\caption{Solution at $t=\SI{50}{\milli\second}$.}
		\label{fig:end_val_fibrosis}
	\end{subfigure}\hfill%
	\caption{Left atrium with fibrosis. The fibrosis pattern (panel B), has been obtained by clipping at level 0.5 the random field (panel A). All elements within the red region in panel (B) are considered ``fibrotic'' and assigned a transverse conductivity of 0 (no transverse conduction). We also report the initial value (panel C) and solution at $t=\SI{50}{\milli\second}$ (panel D) for the CRN ionic model.}
	\label{fig:fibrosis}
\end{figure}

Finally, we repeat the numerical experiments from \cref{sec:fastl_LA} using the same computational mesh 
but in the presence of fibrosis in the cardiac tissue---a common condition in patients with history of atrial fibrillation. Mathematically, fibrosis can be modelled by locally reducing the electric conductivity tensor and altering some ionic currents in the membrane model.
%Indeed fibrosis refers to the scarring of myocardial tissue, which replaces healthy muscle tissue with scar tissue and thus disrupts the normal propagation of electric potential through the heart. 
Here, we incorporate fibrosis only in the conductivity tensor by decoupling cells in the cross-fiber direction~\cite{Gharaviri2020epi,Gharaviri2021PVI,GanPezGhaKraPerSah22}. The spatial heterogeneity is obtained from a spatially-correlated random field $u:\Omega\rightarrow [0,1]$, where we set $\sigma_t=0$ if $u(\bm x)>q$, where $q$ is some quantile value; in our simulations, we use the median. Following \cite[Section 2]{PezQuaPot19}, we generate the random field $u$, shown in \cref{fig:fibrosis_field}, by solving a stochastic differential equation. The corresponding fibrosis pattern is shown in \cref{fig:fibrosis_pattern}.

\begin{figure}	
	\begin{subfigure}[t]{\subfigsizet\textwidth}
		\centering
		\includegraphics[scale=\plotimscalet,trim={0 3mm 0 0},clip]{\currfiledir 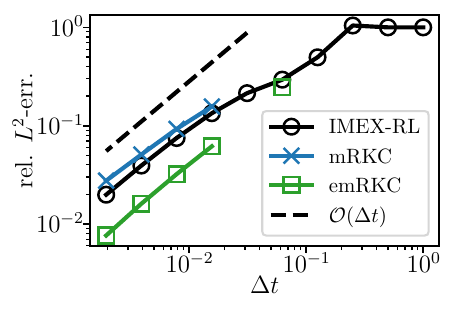}
		\caption{HH model.}
		\label{fig:03_fastl_LA_conv_ref_2_HH_fibrosis}
	\end{subfigure}\hfill%
	\begin{subfigure}[t]{\subfigsizet\textwidth}
		\centering
		\includegraphics[scale=\plotimscalet,trim={0 3mm 0 0},clip]{\currfiledir 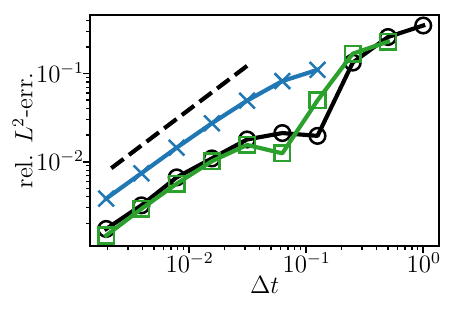}
		\caption{CRN model.}
		\label{fig:03_fastl_LA_conv_ref_2_CRN_fibrosis}
	\end{subfigure}\hfill%
		\begin{subfigure}[t]{\subfigsizet\textwidth}
		\centering
		\includegraphics[scale=\plotimscalet,trim={0 3mm 0 0},clip]{\currfiledir 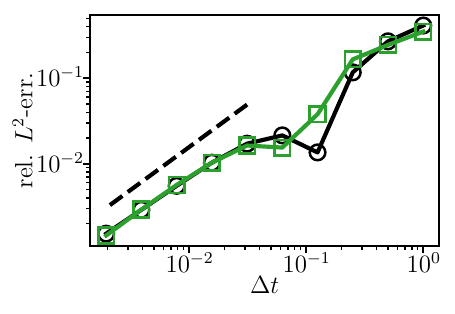}
		\caption{TTP model.}
		\label{fig:03_fastl_LA_conv_ref_2_TTP_fibrosis}
	\end{subfigure}
	\caption{Left atrium with fibrosis. Convergence experiments with different ionic models.}
	\label{fig:03_fastl_LA_conv_fibrosis}
\end{figure}

\begin{figure}	
	\begin{subfigure}[t]{\subfigsizet\textwidth}
		\centering
		\includegraphics[scale=\plotimscalet,trim={0 3mm 0 0},clip]{\currfiledir 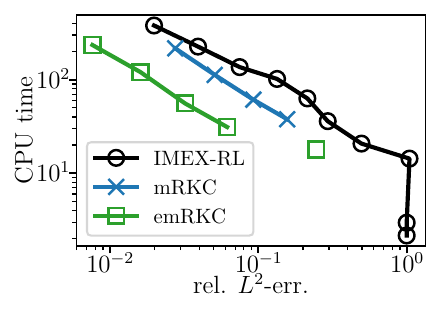}
		\caption{HH model.}
		\label{fig:03_fastl_LA_eff_ref_2_HH_fibrosis}
	\end{subfigure}\hfill%
	\begin{subfigure}[t]{\subfigsizet\textwidth}
		\centering
		\includegraphics[scale=\plotimscalet,trim={0 3mm 0 0},clip]{\currfiledir 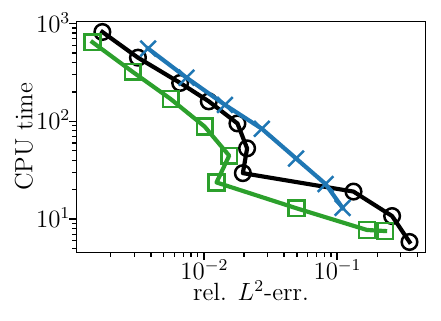}
		\caption{CRN model.}
		\label{fig:03_fastl_LA_eff_ref_2_CRN_fibrosis}
	\end{subfigure}\hfill%
	\begin{subfigure}[t]{\subfigsizet\textwidth}
		\centering
		\includegraphics[scale=\plotimscalet,trim={0 3mm 0 0},clip]{\currfiledir 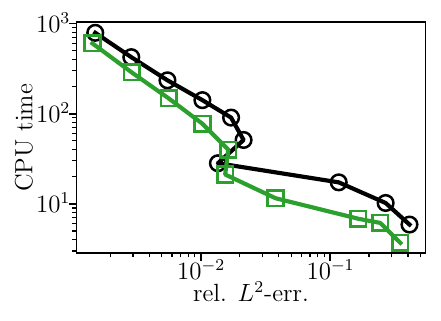}
		\caption{TTP model.}
		\label{fig:03_fastl_LA_eff_ref_2_TTP_fibrosis}
	\end{subfigure}
	\caption{Left atrium with fibrosis. Efficiency experiments with different ionic models.}
	\label{fig:03_fastl_LA_eff_fibrosis}
\end{figure}

The initial value is generated exactly as in \cref{sec:fastl_LA} and is displayed in \cref{fig:init_val_fibrosis}, whereas the solution at the final time $T=\SI{50}{\milli\second}$ is shown in \cref{fig:end_val_fibrosis}. Again in \cref{fig:03_fastl_LA_conv_fibrosis}, we verify first-order convergence for all three methods and models and also observe instabilities in mRKC but not in emRKC. The number of stages taken by the methods is omitted here since it essentially coincides with that in \cref{fig:03_fastl_LA_stages}. In \cref{fig:03_fastl_LA_eff_fibrosis}, we display the results from the efficiency experiments where we again observe that emRKC is slightly faster than IMEX-RL.

\subsection{Parallel performance and scalability}\label{sec:scalability}
To illustrate the inherent parallelism of the explicit emRKC method, we now
conduct a strong scalability experiment where for a fixed problem size we progressively
increase the number of processors. Again we consider the left atrium geometry from \cref{sec:fastl_LA}  and use the CRN ionic model with 21 state variables (20 for the ionic model and 1 for the potential).
Since the FE mesh has about $1.7\times 10^6$ degrees of freedom (dof's), the resulting total problem size is about $36\times 10^6$ dof's.

\begin{figure}
    \centering
    \includegraphics[scale=0.75,trim={0 3mm 0 0},clip]{\currfiledir 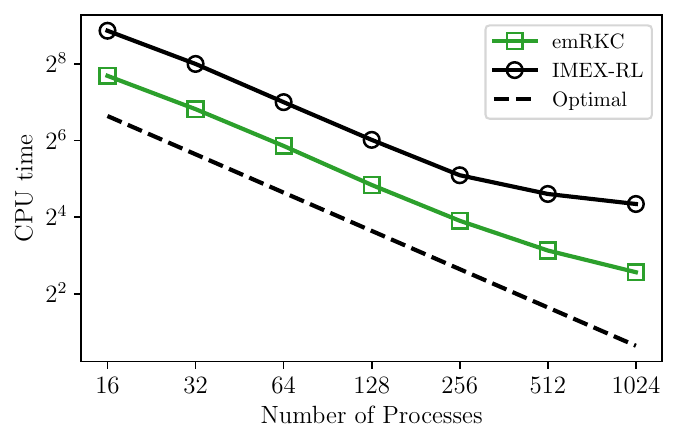}
    \caption{Strong scalability. CPU time vs. the number of processors.}
    \label{fig:scalability}
\end{figure}

In \cref{fig:scalability}, we display for both the emRKC and IMEX-RL methods
the total execution CPU time against the number of processors ranging from 16 to 1024; hence, we start with about $2.2\times 10^6$ dof's 
per processor and end up with only about $35'000$ dof's per processor.
Although both methods yield an optimal speed-up up to 256 processors for this problem size, the parallel efficiency of the implicit IMEX-RL then declines.
Although those preliminary computational results are limited to 1024 processors and a relatively modest problem size, they already illustrate not only that the explicit emRKC method is faster than the IMEX-RL methods, as expected from \cref{fig:prof_03_fastl_emRKC} and \cref{fig:prof_03_fastl_IMEXRL}, but also that the explicit emRKC method exhibits better parallel scalability than the IMEX-RL method.

\section{Conclusion}\label{sec:conclusion}

Explicit stabilized methods, such as RKC, strike a judicious balance between explicit and implicit methods in particular for the time integration of large-scale nonlinear parabolic PDE's. 
The monodomain equation from cardiac electrophysiology is in principle well-suited for explicit stabilized methods due to its moderately stiff diffusive part. However, in the presence of coupled stiff nonlinear
ionic models, which are inherently multiscale in nature, conventional RKC methods can become ineffective.
By adapting the number of stages locally to individual components depending on their respective
degree of stiffness, multirate RKC methods (mRKC) \cite{AGR22} significantly mitigate this efficiency loss, yet without exploiting the particular diagonal structure of ionic models as in the standard Rush--Larsen approach.

To further improve upon the mRKC approach, we have proposed the emRKC method
tailored specifically to the monodomain model. Drawing inspiration from the Rush--Larsen approach, the emRKC method leverages the special diagonal structure of ionic models for which it applies exponential integration. The combination of mRKC with exponential time integration results in a particularly efficient approach for solving the monodomain equations. Moreover, since it is fully explicit, the emRKC is particularly well-suited for implementation on GPU architectures in a parallel HPC environment. 
As existing fully explicit solvers for the monodomain model are always limited to structured grids~\cite{KPD12}, the emRKC thus removes this limitation and extends explicit time integration to unstructured grids, which are crucial for electro-mechanical~\cite{Favino2016impact,quarteroni2022modeling} and sub-cellular~\cite{BidoBEM,deSouzaGaps23,Huynh2023} models.

The stability and accuracy properties of the emRKC method have been analyzed in \cref{sec:stab_analysis,sec:acc_analysis}. In particular, in \cref{thm:stab_emRKC} we show that for a linear scalar test equation the method is stable without constraint on the step size for sufficiently many stages. We also prove in \cref{thm:convergence} that the method is first-order accurate. Both stability and accuracy are also verified numerically. In particular, we show in \cref{sec:eff_stab} through a realistic nonlinear example that emRKC can take up to 1500 times larger time steps than a standard explicit method.

In our numerical experiments, emRKC outperformed or was on par with a standard baseline method, 
which combines IMEX integration for the parabolic component with exponential integration for the ionic model, as such a fair benchmark for comparison. All numerical experiments unequivocally showcased the good performance of emRKC, even in the presence of severe stiffness and nonlinearity. 

The results of our study provide compelling evidence that explicit methods can indeed be exploited effectively for the solution of large-scale, stiff, nonlinear equations and can be a robust and efficient choice for simulating multiscale dynamics. 
Furthermore, the outcomes of this study provide a foundation for the design and implementation of explicit parallel-in-time methods specifically tailored for addressing the monodomain equation.

\appendix

\section{Nonlinear power iteration}\label{app:pow_it}
Here we provide a pseudo-code for the nonlinear power iteration \cite{Lin72,Lin73,Sha91,shampine1980lipschitz,Ver80} used in the RKC, mRKC, and emRKC methods to estimate the spectral radii of the right-hand sides $f$, $f_F$, $f_S$. 
%The method reported here takes inspiration from those proposed in .
The function \textsc{Compute\_Rho} below takes as input parameters the time $t$ and state $y$, for which the spectral radius of $\tfrac{\partial f}{\partial y}(t,y)$ needs to be estimated, the nonlinear function $f$, and an initial guess $v,\rho$ for the dominant eigenvector of $\tfrac{\partial f}{\partial y}(t,y)$ and its spectral radius. If no initial guess $v,\rho$ is available, we simply choose an arbitrary random vector and set $\rho=0$.
The computed values $v,\rho$ returned by \textsc{Compute\_Rho} provide an initial guess for the next function call. 

\IfStandalone{}{\begin{algorithm}}
	\begin{algorithmic}[1]
		\Function{compute\_rho}{$t$, $y$, $f$, $v$, $\rho$}
		\State $n_y=\Vert y\Vert$, $n_v=\Vert v\Vert$
		\State $\epsilon=10^{-8}$, $\delta = n_y\epsilon$, $q=\delta/n_v$
		\State $tol=10^{-2}$, $\rho_\text{old}=0$
		\While{$|\rho-\rho_\text{old}|\geq tol\cdot \rho$}		
            \State $\rho_\text{old}=\rho$, $v_\text{old}=v$
		\State $z=y+q\cdot v_\text{old}$
		\Comment{$\Vert z-y\Vert=\delta$}
		\State $v=f(t,z)-f(t,y)$ \label{algo:powiter_v}
            \Comment{$v\approx \tfrac{\partial f}{\partial y}(t,y)(z-y)=q\tfrac{\partial f}{\partial y}(t,y)v_\text{old}$}
		\State $n_v=\Vert v\Vert$, $q=\delta/n_v$		
		\State $\rho = n_v/\delta$\label{algo:powiter_rho}
		\Comment{$n_v/\delta=\Vert f(t,z)-f(t,y)\Vert/\Vert z-y\Vert$  }
		\EndWhile
		\State \Return $v,\rho$
		\EndFunction
	\end{algorithmic}
	\IfStandalone{}{
		\caption{Nonlinear power iteration}
		\label{algo:powiter}
	\end{algorithm}
}

The function \textsc{Compute\_Rho} indeed corresponds to a power iteration, as \cref{algo:powiter_v} of \cref{algo:powiter} implies that
\begin{equation}
    v\approx \tfrac{\partial f}{\partial y}(t,y)(z-y)=q\tfrac{\partial f}{\partial y}(t,y)v_\text{old},
\end{equation}
while \cref{algo:powiter_rho}
\begin{equation}
    \rho=\frac{n_v}{\delta}=\frac{\Vert f(t,z)-f(t,y)\Vert}{\Vert z-y\Vert}\geq \frac{\langle f(t,z)-f(t,y),z-y \rangle}{\Vert z-y\Vert^2}\approx \frac{\langle \tfrac{\partial f}{\partial y}(t,y) q\cdot v_\text{old},q\cdot v_\text{old} \rangle}{\Vert q\cdot v_\text{old}\Vert^2} = \frac{\langle \tfrac{\partial f}{\partial y}(t,y) v_\text{old},v_\text{old} \rangle}{\Vert v_\text{old}\Vert^2}
\end{equation}
approximately bounds the Rayleigh quotient.
In practice, $\rho$ is subsequently multiplied by a safety factor $\alpha > 1$, typically $\alpha=1.05$.

The added computational cost from \cref{algo:powiter} during time integration scheme is marginal. On the one hand, there is no need to update $v,\rho$ at every time step, every so often, or even just once, is usually sufficient. On the other hand, subsequent calls typically converge very fast ($1-2$ iterations) thanks to the previous values of $v,\rho$ now available as initial guesses.

\ifstandalone
\bibliographystyle{abbrv}
\bibliography{../library}
\fi

\begin{acknowledgement}
\textbf{Acknowledgements} 
The first author thanks the developers of FEniCSx and the Swiss National Supercomputing Centre staff for their support.
\end{acknowledgement}

%\appendix
%\import{src/}{app_ES}

%%-----------------------------
%%      your bibliography
%%-----------------------------
\bibliographystyle{abbrv}
\bibliography{library}
\end{document}